\documentclass{article}
\usepackage{latexsym,amsthm,amssymb,amscd,amsmath}
\newcounter{alphthm}
\newtheorem{thm}{Theorem}

\newtheorem{defn}{Definition}
\newtheorem{prop}{Proposition}
\newtheorem{cor}{Corollary}
\newtheorem{lem}{Lemma}

\newcommand{\be}{\begin{equation}}
\newcommand{\ee}{\end{equation}}
\newcommand{\ben}{\begin{enumerate}}
\newcommand{\een}{\end{enumerate}}
\newcommand{\pa}{{\partial}}

\newcommand{\g}{{\bf g}}
\newcommand{\pxi}{{\pa \over \pa x^i}}
\newcommand{\D}{\delta^d}

\def\beq{\begin{equation}}
\def\eeq{\end{equation}}

\title{On Doubly Warped Product Finsler Manifolds}
\author{E. Peyghan and A. Tayebi }
\begin{document}

\maketitle
\begin{abstract}
In this paper, we introduce horizontal and vertical warped product Finsler manifold. We prove that every C-reducible or proper Berwaldian doubly warped product Finsler manifold is Riemannian.   Then, we find the relation between Riemmanian curvatures of doubly warped product Finsler manifold and  its components, and  consider the cases that this manifold is flat or has scalar flag curvature. We define the doubly warped Sasaki-Matsumoto metric  for warped product manifolds and find a condition under which the horizontal and vertical tangent bundles are totally geodesic. Also,  we obtain some conditions under which a foliated manifold  reduces to  a Reinhart manifold. Finally, we study an almost complex structure on the slit tangent bundle of a doubly warped product Finsler manifold.\\\\
{\bf {Keywords}}:  Doubly warped product manifold, Sasaki-Matsumoto lift metric, Vaisman connection, Reinhart manifold, K\"{a}hler structure.\footnote{ 2000 Mathematics subject Classification: 53C60, 53C25.}
\end{abstract}

\section{Introduction.}
In the Riemannian or semi-Riemannian cases, the doubly warped product of Riemannian (semi-Riemannian) manifolds was studied by many authors \cite{A, BEP, BP, G, U}, and several application to theoretical physics  were given. For instance in \cite{BP}, Beem-Powell considered this product for Lorentzian manifolds. Also, Allison considered causality and global hyperbolicity of doubly warped product and null pseudo convexity of Lorentizian doubly warped product in \cite{A}.

This construction can be extended for Finslerian metrics with some minor restriction. In \cite{As, As1}, Asanov  gave the generalization of the Schwarzschild metric in the Finslerian setting and  obtain some models of relativity theory described through the warped product of Finsler metrics. Then,  Shen  used a construction of warped of Riemannian metrics at the vertical bundle, and obtained a Finslerian warped product metric \cite{S}. Recently,  Kozma-Peter-Varga used the Finsler fundamental functions to define their warped product \cite{Koz}. Then they studied the relationships between the Cartan connection of the doubly warped product manifold  and components of it.

Let $(M_1,F_1)$ and $(M_2,F_2)$ be two Finsler manifolds with $dim M_1=n_1$ and $dim M_2=n_2$ and $f_1: M_1\rightarrow \mathbb{R}^+$ and $f_2: M_2\rightarrow \mathbb{R}^+$ be two smooth functions. Let $\pi_1: M_1\times M_2\rightarrow M_1$ and $\pi_2: M_1\times M_2\rightarrow M_2$ be the natural projection maps. The product manifold $M_1\times M_2$ endowed with the metric $F: TM^\circ_1\times TM^\circ_2\rightarrow\mathbb{R}$ is considered,
\begin{equation}
F(v_1, v_2)=\sqrt{f_2^2(\pi_2(v_2))F_1^2(v_1)+f_1^2(\pi_1(v_1))F_2^2(v_2)},\label{Finsler}
\end{equation}
where $TM^\circ_1=TM_1-\{0\}$ and $TM^\circ_2=TM_2-\{0\}$. The metric defined above is a Finsler metric. The product manifold $M_1\times M_2$ with the metric $F(v)=F(v_1,v_2)$ for $(v_1,v_2)\in TM_1^\circ\times TM_2^\circ$ defined above will be called the doubly warped product of the manifolds $M_1$ and $M_2$ and $f_1$ and $f_2$ will be called the warping functions. We denote this doubly warped by ${}_{f_2}M_1\times_{f_1}M_2$. If either $f_1=1$ or $f_2=1$, but not both, then ${}_{f_2}M_1\times_{f_1}M_2$ becomes a warped product of Finsler manifolds $M_1$ and $M_2$. If both $f_1$ and $f_2$, then we have a product manifold. If neither $f_1$ nor $f_2$ is constant, then we have
a nontrivial (proper) doubly warped product manifold

This paper is arranged as follows: In section 2, we give some  of basic concepts related to Finsler manifolds. In section 3, we introduce the horizontal and vertical distributions on tangent bundle of a doubly warped product Finsler manifold and  construct the Finsler connection on this manifold. Then,  we prove that every C-reducible or proper Berwaldian doubly warped product Finsler manifold reduces to a Riemannian manifold. In section 4, for very two Finsler manifolds $(M_1,F_1)$ and $(M_2,F_2)$,  we introduce the Riemmanian curvature of doubly warped product Finsler manifold $({}_{f_2}M_1\times{}_{f_1}M_2,F)$ and find the relation between it and Riemmanian curvatures of its components $(M_1,F_1)$ and $(M_2,F_2)$. In the cases that $({}_{f_2}M_1\times{}_{f_1}M_2,F)$ is flat or it has the scalar flag curvature, we obtain some results on its components. In section 5, the doubly warped Sasaki-Matsumoto metric ${\bf G}$ is introduced for the doubly warped product Finsler manifold. Then by using the Levi-Civita connection of this metric, we find some conditions under which  $HTM^\circ$ and $VTM^\circ$ are totally geodesic. In section 6, we obtain the Vaisman connection of Riemmanian foliated manifold $(TM^\circ, \mathcal{F}_V,
\textbf{G})$ and  show that it is a Reinhart space if and only if $(M_1,F_1)$ and $(M_2,F_2)$ are Riemannian manifolds. Finally, we define an almost complex structure on the slit tangent bundle of a doubly warped product Finsler manifold and show that this structure  with the doubly warped Sasaki-Matsumoto metric construct an almost Hermitian structure. Then, we prove that $(TM^\circ,\textbf{G}, \textbf{J})$ is a K\"{a}hlerian manifold if and only if the doubly warped horizontal distribution $HTM^\circ$ is integrable.

%-----------------------------------------------------------------------------------------------------------------------------
\section{Preliminary}
%-----------------------------------------------------------------------------------------------------------------------------
Let $M$ be a n-dimensional $ C^\infty$ manifold. Denote by $T_x M $ the tangent space at $x \in M$,  by $TM=\cup _{x \in M} T_x M $ the tangent bundle of $M$, and by $TM^\circ = TM-\{ 0 \}$ the slit tangent bundle on $M$. A {\it Finsler metric} on $M$ is a function $ F:TM \rightarrow [0,\infty)$ which has the following properties: (i) $F$ is $C^\infty$ on $TM^\circ$; (ii) $F$ is positively 1-homogeneous on the fibers of tangent bundle $TM$; (iii) for each $y\in T_xM$, the following quadratic form $\textbf{g}_y$ on $T_xM$  is positive definite,
\[
\textbf{g}_{y}(u,v):={1 \over 2}\frac{\partial^2}{\partial s \partial t} \left[  F^2 (y+su+tv)\right]|_{s,t=0}, \ \
u,v\in T_xM.
\]
Define ${\bf C}_y:T_xM\otimes T_xM\otimes T_xM\rightarrow \mathbb{R}$ by
\[
{\bf C}_{y}(u,v,w):={1 \over 2} \frac{d}{dt}\left[{\bf g}_{y+tw}(u,v)
\right]|_{t=0}, \ \ u,v,w\in T_xM.
\]
The family ${\bf C}:=\{{\bf C}_y\}_{y\in TM^\circ}$  is called the Cartan torsion. It is well known that ${\bf{C}}=0$ if and only if $F$ is Riemannian \cite{BCS}\cite{ShDiff}. For $y\in T_x M_0$, define  mean Cartan torsion ${\bf I}_y$ by ${\bf I}_y(u):=I_i(y)u^i$, where $I_i:=g^{jk}C_{ijk}$, $g^{jk}$ is the inverse of $g_{jk}$ and $u=u^i\frac{\partial}{\partial x^i}|_x$. By Deicke's  Theorem, $F$ is Riemannian  if and only if ${\bf I}_y=0$.

\bigskip

Let $(M, F)$ be a Finsler manifold. Then for  $y \in T_xM_0$, define the  Matsumoto torsion ${\bf M}_y:T_xM\otimes T_xM \otimes T_xM \rightarrow \mathbb{R}$ by ${\bf M}_y(u,v,w):=M_{ijk}(y)u^iv^jw^k$ where
\[
M_{ijk}:=C_{ijk} - {1\over n+1}  \{ I_i h_{jk} + I_j h_{ik} + I_k h_{ij} \},\label{Matsumoto}
\]
$h_{ij}:=FF_{y^iy^j}=g_{ij}-\frac{1}{F^2}g_{ip}y^pg_{jq}y^q$ is the angular metric. A Finsler metric $F$ is said to be C-reducible if ${\bf M}_y=0$ \cite{Mat}.  Matsumoto proves that every Randers metric satisfies that ${\bf M}_y=0$. Later on, Matsumoto-H\={o}j\={o} proves that the converse is true too.  It is remarkable that, a Randers metric $F=\alpha+\beta$ on a manifold $M$ is just a Riemannian metric $\alpha$ perturbated by a one form $\beta$ on $M$ \cite{TP}.

\bigskip

For a Finsler manifold $(M,F)$,  a global vector field ${\bf G}$ is induced by $F$ on $TM_0$, which in a standard coordinate $(x^i,y^i)$ for $TM_0$ is given by ${\bf G}=y^i {{\partial} \over {\partial x^i}}-2G^i(x,y){{\partial} \over {\partial y^i}}$, where \[
G^i:=\frac{1}{4}g^{il}\Big\{\frac{\partial^2[F^2]}{\partial x^k \partial y^l}y^k-\frac{\partial[F^2]}{\partial x^l}\Big\},\ \ y\in T_xM.
\]
The {\bf G} is called the  spray associated  to $(M,F)$.  Then we can define ${\bf B}_y:T_xM\otimes T_xM \otimes T_xM\rightarrow T_xM$  by ${\bf B}_y(u, v, w):=B^i_{\ jkl}(y)u^jv^kw^l{{\partial } \over {\partial x^i}}|_x$ where
\[
B^i_{\ jkl}:={{\partial^3 G^i} \over {\partial y^j \partial y^k \partial y^l}}.
\]
The $\bf B$ is called the Berwald curvature. $F$ is called a Berwald metric  if $\bf{B}=0$.

\bigskip

The Riemann curvature ${\bf R}_y= R^i_{\ k}  dx^k \otimes \pxi|_x :
T_xM \to T_xM$ is a family of linear maps on tangent spaces, defined
by
\begin{equation}\label{TP4}
R^i_{\ k} = 2 {\pa G^i\over \pa x^k}-y^j{\pa^2 G^i\over \pa
x^j\pa y^k} +2G^j {\pa^2 G^i \over \pa y^j \pa y^k} - {\pa G^i \over
\pa y^j} {\pa G^j \over \pa y^k}.
\end{equation}
For a flag $P={\rm span}\{y, u\} \subset T_xM$ with flagpole $y$, the  flag curvature ${\bf K}={\bf K}(P, y)$ is defined by
\begin{equation}\label{TP5}
{\bf K}(P, y):= {\g_y (u, {\bf R}_y(u)) \over \g_y(y, y) \g_y(u,u)
-\g_y(y, u)^2 }.
\end{equation}
We say that a Finsler metric $F$ is   of scalar curvature if for any $y\in T_xM$, the flag curvature ${\bf K}= {\bf K}(x, y)$ is a scalar function on the slit tangent bundle $TM_0$. If ${\bf K}=constant$, then $F$ is said to be of  constant flag curvature.

%-----------------------------------------------------------------------------------------------------------------------------
\section{Doubly Warped Nonlinear Connection}
%-----------------------------------------------------------------------------------------------------------------------------

Let $(M_1,F_1)$ and $(M_2,F_2)$ be two Finsler manifolds. Then the functions
\begin{equation}
(i)\ g_{ij}(x,y)=\frac{1}{2}\frac{\partial^2F_1^2(x,y)}{\partial y^i\partial y^j},\ \ \ (ii)\ g_{\alpha\beta}(u,v)=\frac{1}{2}\frac{\partial^2F_2^2(u,v)}
{\partial v^\alpha\partial v^\beta},\label{metr}
\end{equation}
define a Finsler tensor field of type $(0,2)$ on $TM^\circ_1$ and $TM^\circ_2$, respectively.
 Now let $({}_{f_2}M_1\times{}_{f_1}M_2,F)$ be a doubly warped Finsler manifold and let $\textbf{x}\in M$ and $\textbf{y}\in T_\textbf{x}M$,
 where $\textbf{x}=(x,u)$, $\textbf{y}=(y,v)$, $M=M_1\times M_2$ and $T_\textbf{x}M=T_xM_1\oplus T_uM_2$. Then by using  (\ref{Finsler}) and (\ref{metr}),  we conclude that
\begin{equation}
\Big(\textbf{g}_{ab}(x,u,y,v)\Big)=\Big(\frac{1}{2}\frac{\partial^2F^2(x,u,y,v)}{\partial \textbf{y}^a\textbf{y}^b}\Big)=\left[
\begin{array}{cc}
f_2^2g_{ij}&0\\
0&f_1^2g_{\alpha\beta}
\end{array}
\right],\label{Mat}
\end{equation}
where $\textbf{y}^a=(y^i,v^\alpha)$, $\textbf{y}^b=(y^j,v^\beta)$, $\textbf{g}_{ij}=f_2^2g_{ij}$, $\textbf{g}_{\alpha\beta}=f_1^2g_{\alpha\beta}$, $\textbf{g}_{i\beta}=\textbf{g}_{\alpha j}=0$, $i, j,\ldots\in\{1,\ldots,n_1\}$, $\alpha, \beta,\ldots\in\{1,\ldots,n_2\}$ and $a, b,\ldots\in\{1,\ldots,n_1+n_2\}$.

\bigskip

Now we consider the the spray coefficients of $F_1$, $F_2$ and $F$ as follows
\begin{eqnarray}
G^i(x,y)\!\!\!\!&=&\!\!\!\!\frac{1}{4}g^{ih}\Big(\frac{\partial^2F_1^2}{\partial y^h\partial x^j}y^j-\frac{\partial F_1^2}{\partial x^h}\Big)(x,y),\label{spray1}\\
G^{\alpha}(u,v)\!\!\!\!&=&\!\!\!\!\frac{1}{4}g^{\alpha\gamma}\Big(\frac{\partial^2F_2^2}{\partial v^\gamma\partial u^\beta}v^\beta
-\frac{\partial F_2^2}{\partial u^\gamma}\Big)(u,v),\label{spray2}\\
{\textbf{G}}^a(\textbf{x},\textbf{y})\!\!\!\!&=&\!\!\!\!\frac{1}{4}{\textbf{g}}^{ab}\Big(\frac{\partial^2F^2}{\partial \textbf{y}^b\partial \textbf{x}^c}\textbf{y}^c
-\frac{\partial F^2}{\partial \textbf{x}^b}\Big)(\textbf{x},\textbf{y}).\label{spray3}
\end{eqnarray}
Taking into account the homogeneity of both $F_1^2$ and $F_2^2$, we can derive from (\ref{spray1}) and (\ref{spray2}) that $G^i$ and $G^\alpha$ are positively homogeneous of degree two respect to $(y^i)$ and $(v^\alpha)$, respectively. Hence, by Euler's theorem about the homogeneous functions, we conclude  that
\[
\frac{\partial G^i}{\partial y^j}y^j=2G^i\ \  \textrm{and}  \ \ \frac{\partial G^\alpha}{\partial v^\beta}v^\beta=2G^\alpha.
\]
By setting $a=i$ in (\ref{spray3}),  we have
\be
\textbf{G}^i(x,u,y,v)=\frac{1}{4}\textbf{g}^{ih}\Big(\frac{\partial^2F^2}{\partial y^h\partial x^j}y^j+\frac{\partial^2F^2}{\partial y^h\partial u^\alpha}v^\alpha-\frac{\partial F^2}{\partial x^h}\Big).\label{G}
\ee
Direct calculations give us
\begin{eqnarray*}
\frac{\partial F^2}{\partial x^h}\!\!\!\!&=&\!\!\!\!\ f_2^2\frac{\partial F_1^2}{\partial x^h}+\frac{\partial f_1^2}{\partial x^h}F_2^2\\
\frac{\partial^2F^2}{\partial y^h\partial x^j}\!\!\!\!&=&\!\!\!\!\ f_2^2\frac{\partial^2F_1^2}{\partial y^h\partial x^j}\\
\frac{\partial^2F^2}{\partial y^h\partial u^\alpha}\!\!\!\!&=&\!\!\!\!\ \frac{\partial f_2^2}{\partial u^\alpha}\frac{\partial F_1^2}{\partial y^h}.
\end{eqnarray*}
Putting these equations together $\textbf{g}^{ih}=\frac{1}{f_2^2}g^{ih}$ in (\ref{G}) and using (\ref{spray1}) imply that
\begin{equation}
\textbf{G}^i(x,u,y,v)=G^i(x,y)+\frac{1}{4f_2^2}g^{ih}\Big(\frac{\partial f_2^2}{\partial u^\alpha}\frac{\partial F_1^2}{\partial y^h}v^\alpha-\frac{\partial f_1^2}{\partial x^h}F_2^2\Big).\label{spray4}
\end{equation}
Similarly, by putting $a=\alpha$ in (\ref{spray3}) and using (\ref{spray2}),  we obtain
\begin{equation}
\textbf{G}^\alpha(x,u,y,v)=G^\alpha(u,v)+\frac{1}{4f_1^2}g^{\alpha\gamma}\Big(\frac{\partial f_1^2}{\partial x^j}\frac{\partial F_2^2}{\partial v^\gamma}y^j-\frac{\partial f_2^2}{\partial u^\gamma}F_1^2\Big).\label{spray5}
\end{equation}
Therefore we have $\textbf{G}^a=(\textbf{G}^i,\textbf{G}^\alpha)$, where $\textbf{G}^a$, $\textbf{G}^i$ and $\textbf{G}^\alpha$ are given by (\ref{spray3}), (\ref{spray4}) and (\ref{spray5}), respectively. Now, we put \begin{equation}
(i)\ \ \textbf{G}^a_b:=\frac{\partial \textbf{G}^a}{\partial \textbf{y}^b},\ \ \ (ii)\ \ G^i_j:=\frac{\partial G^i}{\partial y^j},\ \ \
(iii)\ \ G^\alpha_\beta:=\frac{\partial G^\alpha}{\partial v^\beta}.\label{spray6}
\end{equation}
Then we have
\begin{lem}
The coefficients $\textbf{G}^a_b$ defined by (\ref{spray6}) satisfy in the following
\begin{equation}
\Big(\textbf{G}^a_b(x,u,y,v)\Big)=\left[
\begin{array}{cc}
\textbf{G}^i_j(x,u,y,v)&\textbf{G}^\alpha_j(x,u,y,v)\\
\textbf{G}^i_\beta(x,u,y,v)&\textbf{G}^\alpha_\beta(x,u,y,v)
\end{array}
\right],
\end{equation}
where
\begin{eqnarray}
\textbf{G}^i_j(x,u,y,v)\!\!\!\!&:=&\!\!\!\!\frac{\partial\textbf{G}^i}{\partial y^j}=G^i_j-\frac{1}{4f_2^2}\frac{\partial g^{ih}}{\partial y^j}\frac{\partial f_1^2}{\partial x^h}F_2^2
+\frac{1}{2f_2^2}\frac{\partial f_2^2}{\partial u^\alpha}v^\alpha\delta^i_j,\label{spray7}\\
\textbf{G}^i_\beta(x,u,y,v)\!\!\!\!&:=&\!\!\!\!\frac{\partial\textbf{G}^i}{\partial v^\beta}=\frac{1}{4f_2^2}g^{ih}\Big(\frac{\partial f_2^2}{\partial u^\beta}
\frac{\partial F_1^2}{\partial y^h}-\frac{\partial f_1^2}{\partial x^h}\frac{\partial F_2^2}{\partial v^\beta}\Big),\label{spray8}\\
\textbf{G}^\alpha_j(x,u,y,v)\!\!\!\!&:=&\!\!\!\!\frac{\partial\textbf{G}^\alpha}{\partial y^j}=\frac{1}{4f_1^2}g^{\alpha\gamma}\Big(\frac{\partial f_1^2}{\partial x^j}\frac{\partial F_2^2}{\partial v^\gamma}-\frac{\partial f_2^2}{\partial u^\gamma}\frac{\partial F_1^2}{\partial y^j}\Big),\label{spray9}\\
\textbf{G}^\alpha_\beta(x,u,y,v)\!\!\!\!&:=&\!\!\!\!\ \frac{\partial\textbf{G}^\alpha}{\partial v^\beta}=G^\alpha_\beta-\frac{1}{4f_1^2}\frac{\partial g^{\alpha\gamma}}{\partial v^\beta}\frac{\partial f_2^2}{\partial u^\gamma}F_1^2+\frac{1}{2f_1^2}\frac{\partial f_1^2}{\partial x^j}y^j\delta^\alpha_\beta.\label{spray10}
\end{eqnarray}
\end{lem}
\begin{proof}
By using (\ref{spray4}) and (\ref{spray6}),  we have
\begin{equation}
\frac{\partial\textbf{G}^i}{\partial y^j}=\frac{\partial G^i}{\partial y^j}+\frac{1}{4f_2^2}\Big[\frac{\partial g^{ih}}{\partial y^j}\Big(\frac{\partial f_2^2}{\partial u^\alpha}\frac{\partial F_1^2}{\partial y^h}v^\alpha-\frac{\partial f_1^2}{\partial x^h}F_2^2\Big)+g^{ih}\frac{\partial f_2^2}{\partial u^\alpha}\frac{\partial^2F_1^2}{\partial y^j\partial y^h}v^\alpha\Big].\label{proof1}
\end{equation}
But from the part (i) of (\ref{metr}),  we get $\frac{\partial F_1^2}{\partial y^h}=2g_{hk}y^k$. Hence,  we have
\begin{eqnarray}
\nonumber \frac{\partial g^{ih}}{\partial y^j}\frac{\partial F_1^2}{\partial y^h}\!\!\!\!&:=&\!\!\!\! 2\frac{\partial g^{ih}}{\partial y^j}g_{hk}y^k\\
\!\!\!\!&=&\!\!\!\!\ -2g^{ih}y^k\frac{\partial g_{hk}}{\partial y^j}=0.\label{proof2}
\end{eqnarray}
By plugging  (i) of (\ref{metr}) and (\ref{proof2}) in (\ref{proof1}) and using $g^{ih}g_{hj}=\delta^i_j$,  we have (\ref{spray7}). In a similar way, we can obtain (\ref{spray8})-(\ref{spray10}).
\end{proof}

\bigskip

Now, we are going to consider $VTM^{\circ}$, the kernel of the differential of the projection map
\[
\pi:=(\pi_1, \pi_2):TM_1^{\circ}\oplus TM_2^{\circ}\rightarrow M_1\times M_2,
\]
which is a well-defined subbundle of $TTM^\circ$. Locally, $VTM^{\circ}$ is spanned by the natural vector fields $\{\frac{\partial}{\partial y^1},\ldots, \frac{\partial}{\partial y^{n_1}}, \frac{\partial}{\partial v^1},\ldots, \frac{\partial}{\partial v^{n_2}}\}$ and it is called the \textit{doubly warped vertical distribution} on $TM^\circ$. Then, using the functions introduced  by (\ref{spray7})-(\ref{spray10}), the nonholonomic vector fields are defined as follows
\begin{eqnarray}
\frac{\D}{\D x^i}:\!\!\!\!&=&\!\!\!\!\frac{\partial}{\partial x^i}-\textbf{G}^j_i\frac{\partial}{\partial y^j}-\textbf{G}^{\beta}_i\frac{\partial}{\partial v^{\beta}},\label{dec1}\\
\frac{\D}{\D u^\alpha}:\!\!\!\!&=&\!\!\!\!\frac{\partial}{\partial
u^\alpha}-\textbf{G}^j_\alpha\frac{\partial}{\partial
y^j}-\textbf{G}_\alpha^{\beta}\frac{\partial}{\partial v^{\beta}}.
\end{eqnarray}
This make it possible to construct a complementary vector subbundle $HTM^{\circ}$ to $VTM^{\circ}$ in $TTM^{\circ}$,  which is locally presented as follows
\[
HTM^{\circ}:=span\{\frac{\D}{\D x^1},\ldots,\frac{\D}{\D x^{n_1}}, \frac{\D}{\D u^1},\ldots, \frac{\D}{\D u^{n_2}}\}.
\]
$HTM^{\circ}$ is called the \textit{doubly warped horizontal distribution} on $TM^{\circ}$. Thus the tangent bundle of $TM^\circ$ admits the decomposition
\begin{equation}
TTM^\circ=HTM^\circ\oplus VTM^\circ.\label{dec}
\end{equation}
\begin{prop}
Let $({}_{f_1}M_1\times{}_{f_2}M_2,F)$ be a doubly warped product Finsler manifold. Then $\textbf{G}=(\textbf{G}^a_b)$ is the nonlinear connection on $TM=TM_1\oplus TM_2$. Further, we have
\begin{eqnarray*}
\frac{\partial \textbf{G}^i_j}{\partial y^k}y^k+\frac{\partial \textbf{G}^i_j}{\partial v^\gamma}v^\gamma\!\!\!\!&=&\!\!\!\!\textbf{G}^i_j\\
\frac{\partial \textbf{G}^i_\beta}{\partial y^k}y^k+\frac{\partial \textbf{G}^i_\beta}{\partial v^\gamma}v^\gamma\!\!\!\!&=&\!\!\!\! \textbf{G}^i_\beta\\
\frac{\partial \textbf{G}^\alpha_j}{\partial y^k}y^k+\frac{\partial \textbf{G}^\alpha_j}{\partial v^\gamma}v^\gamma\!\!\!\!&=&\!\!\!\!\textbf{G}^\alpha_j\\
\frac{\partial \textbf{G}^\alpha_\beta}{\partial y^k}y^k+\frac{\partial \textbf{G}^\alpha_\beta}{\partial v^\gamma}v^\gamma\!\!\!\!&=&\!\!\!\! \textbf{G}^\alpha_\beta.\\
\end{eqnarray*}
\end{prop}

\begin{defn}
\emph{Using decomposition (\ref{dec}), the doubly warped vertical morphism $v^d:TTM^\circ\rightarrow VTM^\circ$ is defined by \[
v^d:=\frac{\partial}{\partial y^i}\otimes \delta^dy^i+\frac{\partial}{\partial v^\alpha}\otimes \delta^dv^\alpha,
\]
where}
\begin{equation}
(i)\ \delta^dy^i:=dy^i+\textbf{G}^i_jdx^j+\textbf{G}^i_\beta du^\beta,\ \ \ (ii)\ \delta^dv^\alpha:=dv^\alpha+\textbf{G}^\alpha_jdx^j+\textbf{G}^\alpha_\beta du^\beta.\label{new}
\end{equation}
\end{defn}
For this projective morphism, we have
\[
v^d(\frac{\partial}{\partial y^i})=\frac{\partial}{\partial y^i},\ \ \ v^d(\frac{\partial}{\partial v^\alpha})=\frac{\partial}{\partial v^\alpha},\ \ \ v^d(\frac{\D}{\D x^i})=0,\ \ \ v^d(\frac{\D}{\D u^i})=0.
\]
From the above equations,  we get $(v^d)^2=v^d$ and $\ker(v^d)=HTM^\circ$. This mapping is called the \textit{doubly warped vertical projective}.
\begin{defn}
\emph{Using decomposition (\ref{dec}), the \textit{doubly warped horizontal projective} $h^d:TTM^\circ\rightarrow HTM^\circ$ is defined by $h^d=id-v^d$ or}
\[
h^d:=\frac{\D}{\D x^i}\otimes dx^i+\frac{\D}{\D u^\alpha}\otimes du^\alpha.
\]
\end{defn}

For this projective morphism, we have
\[
h^d(\frac{\D}{\D x^i})=\frac{\D}{\D x^i},\ \ \ h^d(\frac{\D}{\D u^\alpha})=\frac{\D}{\D u^\alpha},\ \ \ h^d(\frac{\partial}{\partial y^i})=0,\ \ \ h^d(\frac{\partial}{\partial v^\alpha})=0.
\]
Thus we result that $(h^d)^2=h^d$ and $\ker(h^d)=VTM^\circ$.
\begin{defn}
\emph{Using decomposition (\ref{dec}), the doubly warped almost tangent structure $J^d: HTM^\circ\rightarrow VTM^\circ$ is defined by
\[
J^d: \frac{\partial}{\partial y^i}\otimes dx^i+\frac{\partial}{\partial v^\alpha}\otimes du^\alpha,
\]
or
\[
J^d(\frac{\D}{\D x^i})=\frac{\partial}{\partial y^i},\ \ \ J^d(\frac{\D}{\D u^\alpha})=\frac{\partial}{\partial v^\alpha}, \ \ \ J^d(\frac{\partial}{\partial y^i})=J^d(\frac{\partial}{\partial v^\alpha})=0.
\]
Thus we result that $J^2=0$ and $\ker J=Im J=VTM^\circ$.}
\end{defn}
Here, we introduce some geometrical objects of doubly warped Finsler manifold. In order to simplify the equations, we rewritten the basis of $HTM^\circ$ and $VTM^\circ$ as follows:
\begin{eqnarray*}
\frac{\D}{\D \textbf{x}^a}\!\!\!\!&=&\!\!\!\!\ \frac{\D}{\D x^i}\delta_a^i+\frac{\D}{\D u^\alpha}\delta_a^\alpha,\\
\frac{\partial}{\partial \textbf{y}^a}\!\!\!\!&=&\!\!\!\!\ \frac{\partial}{\partial y^i}\delta_a^i+\frac{\partial}{\partial v^\alpha}\delta_a^\alpha.
\end{eqnarray*}
It is clear that $TTM^\circ=span\{\frac{\D}{\D \textbf{x}^a},\frac{\partial}{\partial \textbf{y}^a}\}$. The Lie brackets of this basis is given by following
\begin{equation}
[\frac{\D}{\D \textbf{x}^a}, \frac{\D}{\D \textbf{x}^b}]=\textbf{R}^c_{\ ab}\frac{\partial}{\partial \textbf{y}^c},\ \ [\frac{\D}{\D \textbf{x}^a}, \frac{\partial}{\partial \textbf{y}^b}]=\textbf{G}^c_{ab}\frac{\partial}{\partial \textbf{y}^c},\ \ [\frac{\partial}{\partial \textbf{y}^a}, \frac{\partial}{\partial \textbf{y}^b}]=0,
\end{equation}
where
\begin{equation}
(i)\ \ \textbf{R}^c_{\ ab}=\frac{\D \textbf{G}^c_a}{\D \textbf{x}^b}-\frac{\D \textbf{G}^c_b}{\D \textbf{x}^a},\ \ \  (ii)\ \ \textbf{G}^c_{\ ab}=\frac{\partial \textbf{G}^c_a}{\partial \textbf{y}^b}.\label{p14}
\end{equation}
\begin{cor}
Let $({}_{f_2}M_1\times{}_{f_1}M_2,F)$ be a doubly warped product Finsler manifold. Then
\[
\textbf{R}^c_{\ ab}=(\textbf{R}^k_{\ ij}, \textbf{R}^k_{\ i\beta}, \textbf{R}^k_{\ \alpha j}, \textbf{R}^k_{\ \alpha \beta}, \textbf{R}^\gamma_{\ ij}, \textbf{R}^\gamma_{\ i\beta}, \textbf{R}^\gamma_{\ \alpha j}, \textbf{R}^\gamma_{\ \alpha\beta})
\]
where
\[
\textbf{R}^k_{\ ij}:=\frac{\D \textbf{G}^k_i}{\D x^j}-\frac{\D \textbf{G}^k_j}{\D x^i},\ \ \ \textbf{R}^k_{\ i\beta}:=\frac{\D \textbf{G}^k_i}{\D u^\beta}-\frac{\D \textbf{G}^k_\beta}{\D x^i}
\]
\[
\textbf{R}^k_{\ \alpha j}:=\frac{\D \textbf{G}^k_\alpha}{\D x^j}-\frac{\D \textbf{G}^k_j}{\D u^\alpha},\ \ \ \textbf{R}^k_{\ \alpha\beta}:=\frac{\D \textbf{G}^k_\alpha}{\D u^\beta}-\frac{\D \textbf{G}^k_\beta}{\D u^\alpha}
\]
\[
\textbf{R}^\gamma_{\ ij}:=\frac{\D \textbf{G}^\gamma_i}{\D x^j}-\frac{\D \textbf{G}^\gamma_j}{\D x^i},\ \ \ \textbf{R}^\gamma_{\ i\beta}:=\frac{\D \textbf{G}^\gamma_i}{\D u^\beta}-\frac{\D \textbf{G}^\gamma_\beta}{\D x^i}
\]
\[
\textbf{R}^\gamma_{\ \alpha j}:=\frac{\D \textbf{G}^\gamma_\alpha}{\D x^j}-\frac{\D \textbf{G}^\gamma_j}{\D u^\alpha},\ \ \ \textbf{R}^\gamma_{\ \alpha\beta}:=\frac{\D \textbf{G}^\gamma_\alpha}{\D u^\beta}-\frac{\D \textbf{G}^\gamma_\beta}{\D u^\alpha}.
\]
\end{cor}

With a simple calculation, we have the following.

\begin{cor}
Let $({}_{f_2}M_1\times{}_{f_1}M_2,F)$ be a doubly warped product Finsler manifold. Suppose that $\textbf{G}=(\textbf{G}^a_b)$ is the nonlinear connection on $TM$. Then
\[
\textbf{G}^c_{ab}=(\textbf{G}^k_{ij}, \textbf{G}^k_{i\beta}, \textbf{G}^k_{\alpha j}, \textbf{G}^k_{\alpha \beta}, \textbf{G}^\gamma_{ij}, \textbf{G}^\gamma_{i\beta}, \textbf{G}^\gamma_{\alpha j}, \textbf{G}^\gamma_{\alpha\beta})
\]
where
\begin{eqnarray*}
\textbf{G}^k_{ij}\!\!\!\!&=&\!\!\!\!\ \frac{\partial\textbf{G}^k_i}{\partial y^j}=G^k_{ij}-\frac{1}{4f_2^2}\frac{\partial^2g^{kh}}{\partial y^j\partial y^i}\frac{\partial f_1^2}{\partial x^h}F_2^2=\textbf{G}^k_{ji},\\
\textbf{G}^k_{i\beta}\!\!\!\!&=&\!\!\!\!\ \frac{\partial\textbf{G}^k_i}{\partial v^\beta}=-\frac{1}{4f_2^2}\frac{\partial g^{kh}}{\partial y^i}\frac{\partial f_1^2}{\partial x^h}\frac{\partial F_2^2}{\partial v^\beta}+\frac{1}{2f_2^2}\frac{\partial f_2^2}{\partial u^\beta}\delta^k_i=\textbf{G}^k_{\beta i},\\
\textbf{G}^k_{\alpha\beta}\!\!\!\!&=&\!\!\!\!\ \frac{\partial\textbf{G}^k_\alpha}{\partial v^\beta}=-\frac{1}{2f_2^2}g_{\alpha\beta}g^{kh}\frac{\partial f_1^2}{\partial x^h}=\textbf{G}^k_{\beta\alpha},\\
\textbf{G}^\gamma_{ij}\!\!\!\!&=&\!\!\!\!\ \frac{\partial \textbf{G}_i^\gamma}{\partial y^j}=-\frac{1}{2f_1^2}g_{ij}g^{\alpha\gamma}\frac{\partial f_2^2}{\partial u^\alpha}=\textbf{G}^\gamma_{ji},\\
\textbf{G}^\gamma_{i\beta} \!\!\!\!&=&\!\!\!\!\ \frac{\partial
\textbf{G}_i^\gamma}{\partial v^\beta}=-\frac{1}{4f_1^2}\frac{\partial
g^{\alpha\gamma}}{\partial v^\beta}\frac{\partial f_2^2}{\partial u^\alpha}\frac{\partial F_1^2}{\partial
y^i}+\frac{1}{2f_1^2}\frac{\partial f_1^2}{\partial x^i}\delta^\gamma_\beta=\textbf{G}^\gamma_{\beta i},\\
\textbf{G}^\gamma_{\alpha\beta} \!\!\!\!&=&\!\!\!\!\ \frac{\partial\textbf{G}^\gamma_\alpha}{\partial v^\beta}=G^\gamma_{\alpha\beta}-\frac{1}{4f_1^2}\frac{\partial^2g^{\gamma \lambda}}{\partial v^\beta\partial v^\alpha}\frac{\partial f_2^2}{\partial u^\lambda}F_1^2=\textbf{G}^\gamma_{\beta\alpha}.
\end{eqnarray*}
\end{cor}
Apart from $\textbf{G}^c_{ab}$, the functions $\textbf{F}^c_{ab}$ are given by
\begin{equation}
\textbf{F}^c_{ab}=\frac{1}{2}\textbf{g}^{ce}\Big(\frac{\D
\textbf{g}_{ea}}{\D \textbf{x}^b}+\frac{\D \textbf{g}_{eb}}{\D
\textbf{x}^a}-\frac{\D \textbf{g}_{ab}}{\D
\textbf{x}^e}\Big).\label{hor}
\end{equation}
\begin{cor}\label{cor}
Let $({}_{f_2}M_1\times{}_{f_1}M_2,F)$ be a doubly warped product Finsler manifold.  Then
\[
\textbf{F}^c_{ab}=(\textbf{F}^k_{ij}, \textbf{F}^k_{i\beta}, \textbf{F}^k_{\alpha j}, \textbf{F}^k_{\alpha \beta}, \textbf{F}^\gamma_{ij}, \textbf{F}^\gamma_{i\beta}, \textbf{F}^\gamma_{\alpha j}, \textbf{F}^\gamma_{\alpha\beta}),
\]
where
\begin{eqnarray}
\textbf{F}^k_{ij}\!\!\!\!&=&\!\!\!\!\  F^k_{ij}-\frac{1}{2}g^{kh}\Big(M^r_j\frac{\partial
g_{hi}}{\partial y^r}+M^r_i\frac{\partial g_{hj}}{\partial
y^r}-M^r_h\frac{\partial g_{ij}}{\partial y^r}\Big)\label{hor2} \\
\textbf{F}^k_{i\beta}  \!\!\!\!&=&\!\!\!\!\ \frac{1}{2f_2^2}g^{kh}\Big(\frac{\partial f_2^2}{\partial
u^\beta}g_{hi}-f_2^2\textbf{G}^r_\beta\frac{\partial
g_{hi}}{\partial y^r}\Big)=\textbf{F}^k_{\beta i}\label{hor3}\\
\textbf{F}^k_{\alpha\beta} \!\!\!\!&=&\!\!\!\!\ -\frac{1}{2f_2^2}g^{kh}\Big(\frac{\partial
f_1^2}{\partial x^h}g_{\alpha\beta}-f_1^2\textbf{G}^\lambda_h\frac{\partial
g_{\alpha\beta}}{\partial v^\lambda}\Big)\label{hor4}\\
\textbf{F}^\gamma_{ij} \!\!\!\!&=&\!\!\!\!\ -\frac{1}{2f_1^2}g^{\gamma\lambda}\Big(\frac{\partial
f_2^2}{\partial u^\lambda}g_{ij}-f_2^2\textbf{G}^r_\lambda\frac{\partial
g_{ij}}{\partial y^r}\Big)\label{hor5}\\
\textbf{F}^\gamma_{i\beta}  \!\!\!\!&=&\!\!\!\!\ \frac{1}{2f_1^2}g^{\gamma\lambda}\Big(\frac{\partial
f_1^2}{\partial
x^i}g_{\beta\lambda}-f_1^2\textbf{G}^\alpha_i\frac{\partial
g_{\beta\lambda}}{\partial v^\alpha}\Big)=\textbf{F}^\gamma_{\beta i}\label{hor6}\\
\textbf{F}^\gamma_{\alpha\beta}  \!\!\!\!&=&\!\!\!\!\ F^\gamma_{\alpha\beta}-\frac{1}{2}g^{\gamma\lambda}\Big(M^\mu_\beta\frac{\partial
g_{\lambda\alpha}}{\partial v^\mu} +M^\mu_\alpha\frac{\partial
g_{\lambda\beta}}{\partial v^\mu}-M^\mu_\lambda\frac{\partial
g_{\alpha\beta}}{\partial v^\mu}\Big)\label{hor7}
\end{eqnarray}
and $F^k_{ij}=\frac{1}{2}g^{kh}(\frac{\delta g_{hi}}{\delta x^j}+\frac{\delta g_{hj}}{\delta x^i}-\frac{\delta g_{ij}}{\delta x^h})$,
 $F^\gamma_{\alpha\beta}=\frac{1}{2}g^{\gamma\lambda}(\frac{\delta g_{\lambda\alpha}}{\delta u^\beta}+\frac{\delta g_{\lambda\beta}}{\delta u^\alpha}
 -\frac{\delta g_{\alpha\beta}}{\delta u^\lambda})$, $M^r_i=\frac{1}{2f_2^2}\frac{\partial f_2^2}{\partial u^\alpha}v^\alpha\delta^r_i
 -\frac{1}{4f_2^2}\frac{\partial g^{rh}}{\partial y^i}\frac{\partial f_1^2}{\partial x^h}F_2^2$ and $M^\mu_\alpha=\frac{1}{2f_1^2}\frac{\partial f_1^2}{\partial x^r}y^r\delta^\mu_\alpha
 -\frac{1}{4f_1^2}\frac{\partial g^{\mu\lambda}}{\partial v^\alpha}\frac{\partial f_2^2}{\partial u^\lambda}F_1^2$.
\end{cor}
\begin{proof}
By using (\ref{hor}), we obtain
\begin{equation}
\textbf{F}^k_{ij}=\frac{1}{2}g^{kh}\Big(\frac{\D g_{hi}}{\D
x^j}+\frac{\D g_{hj}}{\D x^i}-\frac{\D g_{ij}}{\D
x^h}\Big)\label{hor1}
\end{equation}
Since $g_{ij}$ is a function with respect to $(x,y)$, then by (\ref{spray7}) and (\ref{dec1}) we get
\be
\frac{\D g_{hi}}{\D x^j}=\frac{\partial g_{hi}}{\partial
x^j}-G^r_j\frac{\partial g_{hi}}{\partial y^r}-M^r_j\frac{\partial
g_{hi}}{\partial y^r}.\label{M}
\ee
Interchanging $i$ and $j$ in (\ref{M}),  gives us $\frac{\D
g_{hj}}{\D x^i}$. Again, by interchanging $j$ and $h$ in (\ref{M}), we obtain $\frac{\D g_{ij}}{\D x^h}$. By setting these
equation in (\ref{hor1}),  we get (\ref{hor2}). By  similar calculations, we can prove the another relations.
\end{proof}

\begin{lem}\label{lem}
Let $({}_{f_2}M_1\times{}_{f_1}M_2,F)$ be a doubly warped product Finsler manifold.  Then $\textbf{y}^c\textbf{F}^a_{bc}=\textbf{G}^a_b$, where
$\textbf{F}^a_{bc}$ and $\textbf{G}^a_{b}$ are defined by
(\ref{hor}) and (i) of (\ref{spray6}), respectively.
\end{lem}
\begin{proof}
By using (i) of (\ref{metr}), we get
\begin{equation}
(i)\ \ \frac{\partial g_{ij}}{\partial y^r}=\frac{\partial
g_{ir}}{\partial y^j}=\frac{\partial g_{jr}}{\partial y^i},\ \ \
(ii)\ \ \frac{\partial F_1^2}{\partial y^i}=2g_{ij}y^j. \label{p1}
\end{equation}
Since $g_{ij}$ is 0- positive homogenous, then by using Euler's theorem and the part (i) of (\ref{p1}), we obtain
\begin{equation}
y^r\frac{\partial g_{ij}}{\partial y^r}=y^r\frac{\partial
g_{ir}}{\partial y^j}=y^r\frac{\partial g_{jr}}{\partial y^i}=0.
\label{p2}
\end{equation}
Using (\ref{p2}) and $y^jF^k_{ij}$ in (\ref{hor2}) imply that
\begin{equation}
y^j\textbf{F}^k_{ij}=G^k_i-\frac{1}{2}g^{kh}y^jM^r_j\frac{\partial
g_{hi}}{\partial y^r}.\label{p3}
\end{equation}
Direct calculation gives us
\[
y^jM^r_j\frac{\partial g_{hi}}{\partial
y^r}=\frac{1}{2f_2^2}\frac{\partial f_2^2}{\partial
u^\alpha}v^\alpha y^r\frac{\partial g_{hi}}{\partial
y^r}-\frac{1}{4f_2^2}y^j\frac{\partial g^{rs}}{\partial
y^j}\frac{\partial f_1^2}{\partial x^s}F_2^2\frac{\partial
g_{hi}}{\partial y^r}.
\]
Since $g_{hi}$ and $g^{rs}$ are 0-positive homogenous, then from
the above equation we conclude that
\[
y^jM^r_j\frac{\partial g_{hi}}{\partial y^r}=0.
\]
Therefore from (\ref{p3}), we derive
\begin{equation}
y^j\textbf{F}^k_{ij}=G^k_i.\label{IM1}
\end{equation}
From (\ref{spray8}) and (\ref{hor3}) and using (ii) of (\ref{p1}),
(\ref{p2}) and 2-positive homogenously of $F_2^2$ we get
\begin{equation}
v^\beta\textbf{F}^k_{i\beta}= \frac{1}{2f_2^2}\frac{\partial
f_2^2}{\partial
u^\beta}v^\beta\delta^k_i+\frac{1}{4f_2^2}g^{kh}g^{rs}\frac{\partial
g_{hi}}{\partial y^r}\frac{\partial f_1^2}{\partial
x^s}F_2^2.\label{p4}
\end{equation}
On the other hand, we have
\be
g^{kh}g^{rs}\frac{\partial g_{hi}}{\partial
y^r}=g^{kh}g^{rs}\frac{\partial g_{hr}}{\partial
y^i}=-g^{rs}g_{hr}\frac{\partial g^{kh}}{\partial
y^i}=-\frac{\partial g^{ks}}{\partial y^i}.\label{O}
\ee
Setting (\ref{O}) in (\ref{p4}) implies that
\begin{equation}
v^\beta\textbf{F}^k_{i\beta}= \frac{1}{2f_2^2}\frac{\partial
f_2^2}{\partial
u^\beta}v^\beta\delta^k_i-\frac{1}{4f_2^2}\frac{\partial
g^{ks}}{\partial y^i}\frac{\partial f_1^2}{\partial
x^s}F_2^2=\textbf{G}^k_i-G^k_i.\label{IM2}
\end{equation}
From (\ref{IM1}) and (\ref{IM2}), we get the following
\begin{equation}
\textbf{y}^c\textbf{F}^k_{ic}=y^j\textbf{F}^k_{ij}+v^\beta\textbf{F}^k_{i\beta}=\textbf{G}^k_i.
\end{equation}
Similarly we obtain
\[
v^\beta\textbf{F}^\gamma_{i\beta}=\textbf{G}^\gamma_i, \ \ v^\beta\textbf{F}^k_{\alpha\beta}=\textbf{G}^k_\alpha,\ \ v^\beta\textbf{F}^\gamma_{\alpha \beta}=G^\gamma_\alpha,
\]
\[
y^j\textbf{F}^k_{\alpha j}=0, \ \ y^j\textbf{F}^\gamma_{\alpha j}=\textbf{G}^\gamma_\alpha-G^\gamma_\alpha,\  \ \ \ y^j\textbf{F}^\gamma_{ij}=0
\]
These equations give us
\begin{eqnarray*}
&&\textbf{y}^c\textbf{F}^k_{\alpha
c}=y^j\textbf{F}^k_{\alpha
j}+v^\beta\textbf{F}^k_{\alpha\beta}=\textbf{G}^k_\alpha\\
&&\textbf{y}^c\textbf{F}^\gamma_{ic}=y^j\textbf{F}^\gamma_{ij}+v^\beta\textbf{F}^\gamma_{i\beta}=\textbf{G}^\gamma_i\\
&&\textbf{y}^c\textbf{F}^\gamma_{\alpha c}=y^j\textbf{F}^\gamma_{\alpha
j}+v^\beta\textbf{F}^\gamma_{\alpha\beta}=\textbf{G}^\gamma_\alpha.
\end{eqnarray*}
This completes the proof.
\end{proof}
The local components of doubly warped Cartan tensor field of
Manifold $({}_{f_2}M_1\times{}_{f_1}M_2,F)$ is defined by
\[
\textbf{C}^a_{bc}=\frac{1}{2}\textbf{g}^{ae}\frac{\partial \textbf{g}_{bc}}{\partial \textbf{y}^e}.
\]
From this definition, we conclude the following.
\begin{lem}\label{lem3}
Let $({}_{f_2}M_1\times{}_{f_1}M_2,F)$ be a doubly warped product Finsler manifold. Suppose that $C^k_{ij}$ and $C^\gamma_{\alpha\beta}$ be the local components of Cartan tensor field on $M_1$ and $M_2$,
respectively. Then we have
\[
\textbf{C}^c_{ab}=(\textbf{C}^k_{ij}, \textbf{C}^k_{i\beta},
\textbf{C}^k_{\alpha j}, \textbf{C}^k_{\alpha \beta},
\textbf{C}^\gamma_{ij}, \textbf{C}^\gamma_{i\beta},
\textbf{C}^\gamma_{\alpha j}, \textbf{C}^\gamma_{\alpha\beta}),
\]
where
\[
\textbf{C}^k_{ij}=\frac{1}{2}g^{kh}\frac{\partial g_{ij}}{\partial
y^h}=C^k_{ij},\ \ \
\textbf{C}^\gamma_{\alpha\beta}=\frac{1}{2}g^{\gamma\lambda}\frac{\partial
g_{\alpha\beta}}{\partial v^\lambda}=C^\gamma_{\alpha\beta},
\]
and $\textbf{C}^k_{i\beta}=\textbf{C}^k_{\alpha
j}=\textbf{C}^k_{\alpha
\beta}=\textbf{C}^\gamma_{ij}=\textbf{C}^\gamma_{i\beta}=\textbf{C}^\gamma_{\alpha
j}=0$.
\end{lem}

By using the Lemma \ref{lem3}, we can conclude the following.

\begin{cor}\label{cor4}
Let $({}_{f_2}M_1\times{}_{f_1}M_2,F)$ be a doubly warped product
Finsler manifold. Then $({}_{f_2}M_1\times{}_{f_1}M_2,F)$ is a
Riemannian manifold if and only if $(M_1, F_1)$ and $(M_2, F_2)$
are Riemannian manifolds.
\end{cor}

\bigskip

Now, we are going to consider C-reducible doubly warped product Finsler manifold.

\begin{thm}
Every C-reducible doubly warped product Finsler manifold $({}_{f_1}M_1\times{}_{f_2}M_2,F)$ is a Riemannian manifold.
\end{thm}
\begin{proof}
We define the Matsumoto doubly warped tensor $\textbf{M}_{abc}$ as follows:
\be
\textbf{M}_{abc}=\textbf{C}_{abc}-\frac{1}{n+1}\{\textbf{I}_a\textbf{h}_{bc}+\textbf{I}_b\textbf{h}_{ac}
+\textbf{I}_c\textbf{h}_{ab}\},\label{Mat0}
\ee
where $\textbf{I}_a=\textbf{g}^{bc}\textbf{C}_{abc}$, $\textbf{C}_{abc}=\textbf{g}_{cd}\textbf{C}^d_{ab}$ and $\textbf{h}_{ab}=\textbf{g}_{ab}-\frac{1}{F^2}\textbf{y}_a\textbf{y}_b$ is the angular metric. By attention to (\ref{Mat0}) and the relations $\textbf{C}_{ijk}=f_2^2C_{ijk}$ and $\textbf{C}_{\alpha\beta\gamma}=f_1^2C_{\alpha\beta\gamma}$, we obtain
\be
\textbf{M}_{\alpha jk}=-\frac{1}{n+1}\Big\{f_2^2I_\alpha(g_{jk}-\frac{f_2^2}{F^2}y_jy_k)-\frac{f_1^2f_2^2}{F^2}v_\alpha(I_j y_k+I_ky_j)\Big\}.\label{Mat1}
\ee
Contracting (\ref{Mat1}) with $y^jy^k$ implies that
\be
y^jy^k\textbf{M}_{\alpha jk}=-\frac{f_2^2F_1^2}{(n+1)}(1-\frac{f_2^2F_1^2}{F^2})I_\alpha=-\frac{f_1^2f_2^2F_1^2F_2^2}{(n+1)F^2}I_\alpha.\label{Mat2}
\ee
By assumption $\textbf{M}_{\alpha jk}=0$,  and then $I_\alpha=0$, i.e., $(M_2, F_2)$ is a Riemannian manifold. By  similar calculations,  we  can deduce that $(M_1, F_1)$ is a Riemannian manifold. This completes the proof.
\end{proof}

\begin{thm}
Every  proper doubly warped product Finsler manifold $({}_{f_2}M_1\times{}_{f_1}M_2,F)$ with vanishing Berwald curvature  is a Riemannian manifold.
\end{thm}
\begin{proof}
The coefficients of Berwald curvature of a doubly warped product Finsler manifold $({}_{f_2}M_1\times{}_{f_1}M_2,F)$ are given by following:
\begin{eqnarray}
&&\textbf{B}^k_{ijl}=B^k_{ijl}-\frac{1}{4f_2^2}\frac{\partial^3g^{kh}}{\partial y^i\partial y^j\partial y^l}\frac{\partial f_1^2}{\partial x^h}F_2^2
\\
&&\textbf{B}^k_{i\beta l}=-\frac{1}{4f_2^2}\frac{\partial^2g^{kh}}{\partial y^l\partial y^i}\frac{\partial f_1^2}{\partial x^h}\frac{\partial F_2^2}{\partial v^\beta}
\\
&&\textbf{B}^k_{\alpha\beta l}=-\frac{1}{2f_2^2}g_{\alpha\beta}\frac{\partial g^{kh}}{\partial y^l}\frac{\partial f_1^2}{\partial x^h}
\\
&&\textbf{B}^k_{\alpha\beta\lambda}=-\frac{1}{f_2^2}C_{\alpha\beta\lambda}g^{kh}\frac{\partial f_1^2}{\partial x^h}\label{1}
\\
&&\textbf{B}^k_{i\beta\lambda}=-\frac{1}{2f_2^2}\frac{\partial g^{kh}}{\partial y^i}\frac{\partial f_1^2}{\partial x^h}g_{\beta\lambda}
\\
&&\textbf{B}^\gamma_{\alpha\beta\lambda}=B^\gamma_{\alpha\beta\lambda}-\frac{1}{4f_1^2}
\frac{\partial^3g^{\gamma\nu}}{\partial v^\beta\partial v^\alpha\partial v^\lambda}\frac{\partial f_2^2}{\partial u^\nu}F_1^2\label{4}\\
&&\textbf{B}^\gamma_{i\beta\lambda}=-\frac{1}{4f_1^2}\frac{\partial^2g^{\alpha\gamma}}{\partial v^\beta\partial v^\lambda}\frac{\partial f_2^2}{\partial u^\alpha}\frac{\partial F_1^2}{\partial y^i}\\
&&\textbf{B}^\gamma_{ij\lambda}=-\frac{1}{2f_1^2}g_{ij}\frac{\partial g^{\alpha\gamma}}{\partial v^\lambda}\frac{\partial f_2^2}{\partial u^\alpha},\label{3}\\
&&\textbf{B}^\gamma_{ijk}=-\frac{1}{f_1^2}C_{ijk}g^{\alpha\gamma}\frac{\partial f_2^2}{\partial u^\alpha}\label{2a}\\
&&\textbf{B}^\gamma_{i\beta k}=-\frac{1}{2f_1^2}\frac{\partial g^{\alpha\gamma}}{\partial v^\beta}\frac{\partial f_2^2}{\partial u^\alpha}g_{ik}.
\end{eqnarray}
If $({}_{f_2}M_1\times{}_{f_1}M_2,F)$ is Berwaldian, then we have $\textbf{B}^d_{abc}=0$. By (\ref{1}), we  get
\be
C_{\alpha\beta\lambda}g^{kh}\frac{\partial f_1^2}{\partial x^h}=0.\label{B1}
\ee
Multiplying (\ref{B1}) with $g_{kr}$ implies that
\be
C_{\alpha\beta\lambda}\frac{\partial f_1^2}{\partial x^r}=0.\label{B2}
\ee
By  (\ref{B2}),  if $f_1$ is not constant then we result that $C_{\alpha\beta\lambda}=0$, i.e., $(M_2, F_2)$ is Riemannian. In the similar way, from (\ref{2a}) we conclude that if $f_2$ is non constant then $(M_1, F_1)$ is Riemannian. 
\end{proof}
\begin{thm}
Let $({}_{f_2}M_1\times{}_{f_1}M_2,F)$ be a doubly warped product Finsler manifold and $f_1$ is constant on $M_1$ ($f_2$ is constant on $M_2$). Then $({}_{f_2}M_1\times{}_{f_1}M_2,F)$ is Berwaldian if and only if $M_1$ is Riemannian, $M_2$ is Berwaldian and $\frac{\partial g^{\alpha\gamma}}{\partial v^\lambda}\frac{\partial f_2^2}{\partial u^\alpha}=0$ ($M_2$ is Riemannian, $M_1$ is Berwaldian and $\frac{\partial g^{ij}}{\partial y^k}\frac{\partial f_1^2}{\partial x^i}=0$).
\end{thm}
\begin{proof}
Let $({}_{f_2}M_1\times{}_{f_1}M_2,F)$ be a Berwaldian manifold and $f_1$ is constant on $M_1$. Then  from (\ref{2a}) we result that $C_{ijk}=0$, i.e., $(M_1, F_1)$ is Riemannian. Also, (\ref{3}) gives us $\frac{\partial g^{\alpha\gamma}}{\partial v^\lambda}\frac{\partial f_2^2}{\partial u^\alpha}=0$. Differentiating this equation with respect to $(v^\beta)$ we deduce $\frac{\partial^2 g^{\alpha\gamma}}{\partial v^\lambda\partial v^\beta}\frac{\partial f_2^2}{\partial u^\alpha}=0$ and consequently $\frac{\partial^3 g^{\alpha\gamma}}{\partial v^\lambda\partial v^\beta\partial v^\mu}\frac{\partial f_2^2}{\partial u^\alpha}=0$. Setting this equation in (\ref{4}) we derive $B^\gamma_{\alpha\beta\lambda}=0$, i.e., $(M_2, F_2)$ is Berwaldian. In the similar way, we can prove the converse of this assertion.
\end{proof}
\begin{cor}
Let $(M_1\times{}_{f_1}M_2,F)$ be a proper warped product Finsler manifold. Then $(M_1\times{}_{f_1}M_2,F)$ is Berwaldian if and only if $M_2$ is Riemannian, $M_1$ is Berwaldian and $\frac{\partial g^{ij}}{\partial y^k}\frac{\partial f_1^2}{\partial x^i}=0$.
\end{cor}
%-----------------------------------------------------------------------------------------------------------------------------
\section{Riemannian Curvature of a Doubly Warped Product Manifold}
%-----------------------------------------------------------------------------------------------------------------------------
The Riemannian curvature of a doubly warped product Finsler manifold
$({}_{f_2}M_1\times{}_{f_1}M_2,F)$ with respect to Berwald
connection is given by
\begin{equation}
\textbf{R}^{\ a}_{b\ cd}=\frac{\D \textbf{F}^a_{bc}}{\D \textbf{x}^d}-\frac{\D \textbf{F}^a_{bd}}{\D \textbf{x}^c}+\textbf{F}^a_{de}\textbf{F}^e_{bc}
-\textbf{F}^a_{ce}\textbf{F}^e_{bd}.\label{cur}
\end{equation}
For the definition of  Berwald connection see \cite{TAE} and \cite{TN}.

\begin{lem}
Let $({}_{f_2}M_1\times{}_{f_1}M_2,F)$ be a doubly warped product Finsler manifold. Then
\[
\textbf{R}^a_{\ cd}=\textbf{y}^b\textbf{R}^{\ a}_{b\ cd},
\]
where $\textbf{R}^a_{\ cd}$ and $\textbf{y}^b\textbf{R}^{\ a}_{b\ cd}$ are given by (\ref{p14}) and (\ref{cur}).
\end{lem}
\begin{proof}
By using (\ref{cur}), we have
\begin{equation}
\textbf{y}^b\textbf{R}^{\ i}_{b\ kl}=\textbf{y}^b\frac{\D \textbf{F}^i_{bk}}{\D \textbf{x}^l}-\textbf{y}^b\frac{\D \textbf{F}^i_{bl}}{\D \textbf{x}^k}+\textbf{y}^b\textbf{F}^i_{le}\textbf{F}^e_{bk}.\label{p11}
-\textbf{y}^b\textbf{F}^i_{ke}\textbf{F}^e_{bl}.
\end{equation}
By using Corollary \ref{cor} and Lemma \ref{lem}, we obtain
\begin{equation}
\textbf{y}^b\frac{\D \textbf{F}^i_{bk}}{\D \textbf{x}^l}
=\frac{\D \textbf{G}^i_k}{\D x^l}+\textbf{F}^i_{jk}\textbf{G}^j_l+\textbf{F}^i_{\beta k}\textbf{G}^\beta_l,\ \ \ \textbf{y}^b\textbf{F}^i_{le}\textbf{F}^e_{bk}=\textbf{F}^i_{lh}\textbf{G}^h_k+\textbf{F}^i_{l\gamma}\textbf{G}_k^\gamma.
\label{p12}
\end{equation}
Interchanging $i$ and $j$ in (\ref{p12}) implies that
\begin{equation}
\textbf{y}^b\frac{\D \textbf{F}^i_{bl}}{\D \textbf{x}^k}
=\frac{\D \textbf{G}^i_l}{\D x^k}+\textbf{F}^i_{jl}\textbf{G}^j_k+\textbf{F}^i_{\beta l}\textbf{G}^\beta_k,\ \ \ \textbf{y}^b\textbf{F}^i_{ke}\textbf{F}^e_{bl}=\textbf{F}^i_{kh}\textbf{G}^h_l+\textbf{F}^i_{k\gamma}\textbf{G}_l^\gamma.
\label{p13}
\end{equation}
Plugging  (\ref{p12}) and (\ref{p13}) in (\ref{p11}), give us $\textbf{y}^b\textbf{R}^{\ i}_{b\ kl}=\textbf{R}^i_{\ kl}$. In the similar way, we can obtain this relation for another indices.
\end{proof}

\bigskip

Using (\ref{cur}), we are going to compute the Riemannian curvature of a doubly warped product Finsler manifold.

\begin{lem}
Let $({}_{f_2}M_1\times{}_{f_1}M_2,F)$ be a doubly warped product Finsler manifold. Then  the Riemannian curvature of a doubly warped product Finsler manifold is given by following
\begin{eqnarray}
\textbf{R}^{\ i}_{j\ kl}\!\!\!\!&=&\!\!\!\!R^{\ i}_{j\ kl}-A_{(kl)}\Big\{M^r_l\frac{\partial F^i_{jk}}{\partial y^r}
+\frac{\D M^i_{jk}}{\D x^l}+F^i_{lh}M^h_{jk}+M^i_{lh}F^h_{jk}-M^i_{lh}M^h_{jk}\nonumber\\
\!\!\!\!&&\!\!\!\!+\frac{1}{4f_1^2f_2^2}g^{is}g^{\alpha\gamma}(g_{sl}\frac{\partial f_2^2}{\partial u^\gamma}
-f_2^2\textbf{G}^r_\gamma\frac{\partial g_{sl}}{\partial y^r})(g_{jk}\frac{\partial f_2^2}{\partial u^\alpha}-f_2^2\textbf{G}^h_\alpha\frac{\partial g_{jk}}{\partial y^h})\Big\}\label{cur0}
\end{eqnarray}
\begin{eqnarray}
\textbf{R}^{\ i}_{\alpha\
kl}\!\!\!\!&=&\!\!\!\! A_{(kl)}\Big\{\frac{1}{2f_2^2}\frac{\D}{\D
x^l}(\frac{\partial f_2^2}{\partial
u^\alpha}\delta^i_k-f_2^2\textbf{G}^r_\alpha g^{ih}\frac{\partial
g_{hk}}{\partial
y^r})+\frac{1}{2f_2^2}(F^i_{rl}-M^i_{rl})(\frac{\partial
f_2^2}{\partial u^\alpha}\delta^r_k\nonumber\\
\!\!\!\!&&\!\!\!\!-f_2^2\textbf{G}^m_\alpha g^{rh}\frac{\partial
g_{hk}}{\partial y^m})+\frac{1}{4f_1^2f_2^2}(\frac{\partial
f_2^2}{\partial u^\beta}\delta^i_l-f_2^2\textbf{G}^r_\beta
g^{ih}\frac{\partial g_{hl}}{\partial y^r})(\frac{\partial
f_1^2}{\partial
x^k}\delta^\beta_\alpha\nonumber\\
\!\!\!\!&&\!\!\!\!-f_1^2\textbf{G}^\mu_k
g^{\beta\lambda}\frac{\partial g_{\alpha\lambda}}{\partial
v^\mu})\Big\}
\end{eqnarray}
\begin{eqnarray}
\textbf{R}^{\ i}_{j\
\beta\lambda}\!\!\!\!&=&\!\!\!\!  A_{(\beta \lambda)}\Big\{\frac{\D}{\D
u^\lambda}\Big(\frac{1}{2f_2^2}(\frac{\partial f_2^2}{\partial
u^\beta}\delta^i_j-f_2^2\textbf{G}^r_\beta
g^{ih}\frac{\partial g_{hj}}{\partial
y^r})\Big)+\frac{1}{4(f_2^2)^2}(\frac{\partial
f_2^2}{\partial u^\lambda}\delta^i_r\nonumber\\
\!\!\!\!&&\!\!\!\!-f_2^2\textbf{G}^m_\lambda g^{ih}\frac{\partial
g_{hr}}{\partial y^m})(\frac{\partial f_2^2}{\partial
u^\beta}\delta^r_j-f_2^2\textbf{G}^l_\beta
g^{rs}\frac{\partial g_{sj}}{\partial
y^l})-\frac{1}{4f_1^2f_2^2}(\frac{\partial
f_1^2}{\partial x^h}g^{ih}g_{\alpha\lambda}\nonumber\\
\!\!\!\!&&\!\!\!\!-f_1^2\textbf{G}^\mu_h g^{ih}\frac{\partial
g_{\alpha\lambda}}{\partial v^\mu})(\frac{\partial f_1^2}{\partial
x^j}\delta^\alpha_\beta-f_1^2\textbf{G}^\nu_j
g^{\alpha\gamma}\frac{\partial g_{\beta\gamma}}{\partial
v^\nu})\Big\}
\end{eqnarray}
\begin{eqnarray}
\textbf{R}^{\ i}_{\alpha\ \beta l}\!\!\!\!&=&\!\!\!\!-\frac{\D}{\D
u^\beta}\Big(\frac{1}{2f_2^2}(\frac{\partial f_2^2}{\partial
u^\alpha}\delta^i_l-f_2^2\textbf{G}^r_\alpha g^{ih}\frac{\partial
g_{hl}}{\partial y^r})\Big)+\frac{1}{2f_2^2}\frac{\D}{\D
x^l}(\frac{\partial f_1^2}{\partial
x^h}g^{ih}g_{\alpha\beta}\nonumber\\
\!\!\!\!&&\!\!\!\!-f_1^2\textbf{G}^\lambda_h g^{ih}\frac{\partial
g_{\alpha\beta}}{\partial
v^\lambda})-\frac{1}{4(f_2^2)^2}(\frac{\partial f_2^2}{\partial
u^\beta}\delta^i_r-f_2^2\textbf{G}^m_\beta g^{ih}\frac{\partial
g_{hr}}{\partial y^m})(\frac{\partial f_2^2}{\partial
u^\alpha}\delta^r_l\nonumber\\
\!\!\!\!&&\!\!\!\!-f_2^2\textbf{G}^s_\alpha g^{rn}\frac{\partial
g_{nl}}{\partial y^s})+\frac{1}{4f_1^2f_2^2}(\frac{\partial
f_1^2}{\partial x^h}g^{ih}g_{\mu\beta}-f_1^2\textbf{G}^\lambda_h
g^{ih}\frac{\partial g_{\mu\beta}}{\partial
v^\lambda})(\frac{\partial f_1^2}{\partial
x^l}\delta^\mu_\alpha\nonumber\\
\!\!\!\!&&\!\!\!\!-f_1^2\textbf{G}^\nu_l
g^{\mu\gamma}\frac{\partial g_{\gamma\alpha}}{\partial
v^\nu})-\frac{1}{2f_2^2}g^{rh}(F^i_{rl}-M^i_{rl})(\frac{\partial
f_1^2}{\partial
x^h}g_{\alpha\beta}-f_1^2\textbf{G}^\lambda_h\frac{\partial
g_{\alpha\beta}}{\partial
v^\lambda})\nonumber\\
\!\!\!\!&&\!\!\!\!+\frac{1}{2f_2^2}(\frac{\partial f_2^2}{\partial
u^\mu}\delta^i_l-f_2^2\textbf{G}^r_\mu g^{ih}\frac{\partial
g_{hl}}{\partial y^r})(F^\mu_{\alpha\beta}-M^\mu_{\alpha\beta})
\end{eqnarray}
\begin{eqnarray}
\textbf{R}^{\ \gamma}_{j\
\beta\lambda}\!\!\!\!&=&\!\!\!\!   A_{(\beta \lambda)}\Big\{\frac{1}{2f_1^2}\frac{\D}{\D
u^\lambda}(\frac{\partial f_1^2}{\partial
x^j}\delta^\gamma_\beta-f_1^2\textbf{G}^\alpha_j
g^{\gamma\mu}\frac{\partial g_{\beta\mu}}{\partial
v^\alpha})+\frac{1}{2f_1^2}(F^\gamma_{\alpha\lambda}-M^\gamma_{\alpha\lambda})(\frac{\partial
f_1^2}{\partial x^j}\delta^\alpha_\beta\nonumber\\
\!\!\!\!&&\!\!\!\!-f_1^2\textbf{G}^\nu_j
g^{\alpha\mu}\frac{\partial g_{\beta\mu}}{\partial
v^\nu})+\frac{1}{4f_1^2f_2^2}(\frac{\partial f_1^2}{\partial
x^r}\delta^\gamma_\lambda-f_1^2\textbf{G}^\alpha_r
g^{\gamma\beta}\frac{\partial g_{\lambda\beta}}{\partial
v^\alpha})(\frac{\partial f_2^2}{\partial
u^\beta}\delta^r_j\nonumber\\
\!\!\!\!&&\!\!\!\!-f_2^2\textbf{G}^m_\beta g^{rh}\frac{\partial
g_{hj}}{\partial y^m})\Big\}
\end{eqnarray}
\begin{eqnarray}
\textbf{R}^{\ \gamma}_{\alpha\
kl}\!\!\!\!&=&\!\!\!\!\  A_{(kl)}    \Big\{\frac{\D}{\D
x^l}\Big(\frac{1}{2f_1^2}(\frac{\partial f_1^2}{\partial
x^k}\delta^\gamma_\alpha-f_1^2\textbf{G}^\mu_k
g^{\gamma\lambda}\frac{\partial g_{\alpha\lambda}}{\partial
v^\mu})\Big)+\frac{1}{4(f_1^2)^2}(\frac{\partial
f_1^2}{\partial x^l}\delta^\gamma_\beta\nonumber\\
\!\!\!\!&&\!\!\!\!-f_1^2\textbf{G}^\alpha_l
g^{\gamma\lambda}\frac{\partial g_{\beta\lambda}}{\partial
v^\alpha})(\frac{\partial f_1^2}{\partial
x^k}\delta^\beta_\alpha-f_1^2\textbf{G}^\nu_k
g^{\beta\mu}\frac{\partial g_{\alpha\mu}}{\partial
v^\nu})-\frac{1}{4f_1^2f_2^2}g^{\gamma\lambda}(\frac{\partial
f_2^2}{\partial u^\lambda}g_{lr}\nonumber\\
\!\!\!\!&&\!\!\!\!-f_2^2\textbf{G}^m_\lambda\frac{\partial
g_{lr}}{\partial y^m})(\frac{\partial f_2^2}{\partial
u^\alpha}\delta^r_k-f_2^2\textbf{G}^s_\alpha
g^{rh}\frac{\partial g_{hk}}{\partial
y^s})\Big\}
\end{eqnarray}
\begin{eqnarray}
\textbf{R}^{\ \gamma}_{\alpha\
\beta\lambda}\!\!\!\!&=&\!\!\!\!\ R^{\ \gamma}_{\alpha\
\beta\lambda}-  A_{(\beta \lambda)}\Big\{M^\kappa_\lambda(\frac{\partial
F^\gamma_{\alpha\beta}}{\partial v^\kappa}+\frac{\D M^\gamma_{\alpha\beta}}{\D
u^\lambda})+F^\gamma_{\lambda\mu}M^\mu_{\alpha\beta}
+M^\gamma_{\lambda\mu}(F^\mu_{\alpha\beta}-M^\mu_{\alpha\beta})\nonumber\\
\!\!\!\!&&\!\!\!\!+\frac{1}{4f_1^2f_2^2}g^{\gamma\nu}g^{rs}(g_{\nu\lambda}\frac{\partial f_1^2}{\partial x^s}
-f_1^2\textbf{G}^\kappa_s\frac{\partial g_{\nu\lambda}}{\partial v^\kappa})(g_{\alpha\beta}\frac{\partial f_1^2}{\partial x^r}
-f_1^2\textbf{G}^\mu_r\frac{\partial g_{\alpha\beta}}{\partial v^\mu})\Big\},
\end{eqnarray}
where $M^i_{jk}=\frac{1}{2}g^{ih}(M^r_k\frac{\partial
g_{hj}}{\partial y^r}+M^r_j\frac{\partial g_{hk}}{\partial
y^r}-M^r_h\frac{\partial g_{jk}}{\partial y^r})$,
$M^\gamma_{\alpha\beta}=\frac{1}{2}g^{\gamma
\nu}(M^\mu_\beta\frac{\partial g_{\nu\alpha}}{\partial
v^\mu}+M^\mu_\alpha\frac{\partial g_{\nu\beta}}{\partial
v^\mu}-M^\mu_\nu\frac{\partial g_{\alpha\beta}}{\partial
v^\mu})$ and  $A_{(ij)}$ denotes the interchange of indices $i$, $j$ and subtraction.
\end{lem}
\begin{proof}
By (\ref{cur}),  we have
\begin{equation}
\textbf{R}^{\ i}_{j\ kl}=\frac{\D \textbf{F}^i_{jk}}{\D x^l}-\frac{\D F^i_{jl}}{\D x^k}+F^i_{lh}F^h_{jk}+F^i_{l\gamma}F^\gamma_{jk}-F^i_{kh}F^h_{jl}-F^i_{k\gamma}F^\gamma_{jl}.\label{cur1}
\end{equation}
By using (\ref{hor2}), we get
\[
\textbf{F}^i_{jk}=F^i_{jk}-M^i_{jk}.
\]
Since $F^k_{ij}$ is a function with respect to $(x,y)$, then by  (\ref{spray7}) and
(\ref{dec1}) we derive that
\begin{equation}
\frac{\D \textbf{F}^i_{jk}}{\D x^l}=\frac{\delta F^i_{jk}}{\delta
x^l}-M^r_l\frac{\partial F^i_{jk}}{\partial y^r}-\frac{\D
M^i_{jk}}{\D x^l}.\label{cur2}
\end{equation}
Interchanging $k$ and $l$ in (\ref{cur2}) implies that
\begin{equation}
\frac{\D \textbf{F}^i_{jl}}{\D x^k}=\frac{\delta F^i_{jl}}{\delta
x^k}-M^r_k\frac{\partial F^i_{jl}}{\partial y^r}-\frac{\D
M^i_{jl}}{\D x^k}.\label{cur3}
\end{equation}
By plugging (\ref{hor2}), (\ref{hor3}), (\ref{hor5}), (\ref{cur2})
and (\ref{cur3}) in (\ref{cur1}), we can obtain (\ref{cur0}). In the similar way,
we can obtain this relation for another indices.
\end{proof}

\begin{thm}
Let $({}_{f_2}M_1\times{}_{f_1}M_2,F)$ be a flat doubly warped
product Finsler manifold. Then
\begin{description}
  \item[(i)] if $(M_1, F_1)$ is Riemannian then the components of the
Riemannian Curvature of $M_1$ are as follows:
\begin{equation}
\textbf{R}^{\ i}_{j\ kl}=\frac{||grad
f_2||^2}{f_1^2}(\delta^i_lg_{jk}-\delta^i_kg_{jl}).\label{con}
\end{equation}
  \item[(ii)] if $(M_2, F_2)$ is Riemannian then the components of the
Riemannian Curvature of $M_2$ are as follows:
\begin{equation}
\textbf{R}^{\ \gamma}_{\alpha\ \beta\lambda}=\frac{||grad
f_1||^2}{f_2^2}(\delta^\gamma_\lambda g_{\alpha
\beta}-\delta^\gamma_\beta g_{\alpha\lambda}).
\end{equation}
\end{description}
\end{thm}
\begin{proof}
Since the proof of (ii) similar to (i), then we only prove (i).
Let $(M_1, F_1)$ be a Riemannian manifold. Then $g_{ij}$ is a
function of $(x)$, only. Therefore we have $M^i_{jk}=0$. Also, the
function $F^i_{jk}$ independent of $(y)$. By using
(\ref{cur0}), we conclude that
\[
\textbf{R}^{\ i}_{j\ kl}=R^{\ i}_{j\
kl}-\frac{1}{4f_1^2f_2^2}(\delta^i_lg_{jk}-\delta^i_kg_{jl})g^{\alpha\gamma}\frac{\partial
f_2^2}{\partial u^\gamma}\frac{\partial f_2^2}{\partial u^\alpha}.
\]
But we have
\[
g^{\alpha\gamma}\frac{\partial f_2^2}{\partial u^\gamma}\frac{\partial f_2^2}{\partial u^\alpha}=4f_2^2||grad
f_2||^2.
\]
Hence the above equation rewritten as follows
\begin{equation}
\textbf{R}^{\ i}_{j\ kl}=R^{\ i}_{j\ kl}-\frac{||grad
f_2||^2}{f_1^2}(\delta^i_lg_{jk}-\delta^i_kg_{jl}).\label{con1}
\end{equation}
Since $({}_{f_2}M_1\times{}_{f_1}M_2,F)$ is a flat manifold, then we
have $\textbf{R}^{\ i}_{j\ kl}=0$. Therefore,  (\ref{con1}) gives
us (\ref{con}).
\end{proof}

\bigskip

Now, let $(M_1, F_1)$ is a Riemannian manifold and $f_1$ is a
scalar function on $M_1$. Then from (\ref{con}),  we have
\[
\textbf{R}^{\ i}_{j\ kl}=K_1||gradf_2||^2(\delta^i_lg_{jk}-\delta^i_kg_{jl}),
\]
where $K_1$ is constant. Since $gradf_2$ is independent of $(x)$,
then it is a constant function on $M_1$. The similar argument is hold, if
$(M_2,F_2)$ is Riemannian and $f_2$ is constant on $M_2$.
Therefore we have the following corollary.
\begin{cor}\label{cor5}
Let $({}_{f_2}M_1\times{}_{f_1}M_2,F)$ be a flat doubly warped
product Finsler manifold. Then
\begin{description}
  \item[(i)] if $(M_1, F_1)$ is a Riemannian manifold and $f_1$ is
constant on $M_1$, then $M_1$ is a space of positive constant
curvature $K_1||gradf_2||^2$;
  \item[(ii)] if $(M_2, F_2)$ is a Riemannian manifold and $f_2$ is
constant on $M_2$, then $M_2$ is a space of positive constant
curvature $K_2||gradf_1||^2$.
\end{description}
\end{cor}

\bigskip

By the corollaries \ref{cor4} and  \ref{cor5}, we conclude the following.
\begin{cor}
Let $({}_{f_2}M_1\times{}_{f_1}M_2,F)$ be a flat doubly warped
product Riemannian manifold. Then
\begin{description}
  \item[(i)] if $f_1$ is constant on $M_1$, then $M_1$ has positive constant curvature $K_1||gradf_2||^2$ and $M_2$ is a flat manifold;
  \item[(ii)] if $f_2$ is constant on $M_2$, then $M_2$ has positive constant curvature $K_2||gradf_1||^2$ and $M_1$ is a flat manifold.
\end{description}
\end{cor}

\bigskip

The flag curvature of a Finsler metric which plays the central role in Finsler geometry, is  called a Riemannian quantity  because it is a natural  extension of sectional curvature in Riemannian geometry. For a Finsler manifold $(M, F)$, the flag curvature is  a function ${\bf K}(P, y)$ of tangent planes $P\subset T_xM$ and  directions $y\in P$. The Finsler metric $F$  is said to be  of scalar flag curvature if the flag curvature ${\bf K}(P, y)={\bf K}(x, y)$ is independent of flags $P$ associated with any fixed flagpole $y$ \cite{ShDiff}.

\begin{thm}
Let $(M_1,F_1)$ be a Riemannian \textbf{ma}nifold and
$({}_{f_2}M_1\times{}_{f_1}M_2,F)$ be a doubly warped product
Finsler space of scalar flag curvature $\lambda_1(x,u,y,v)$. Then
$(M_1, F_1)$ has constant curvature $K_1$ if and only if
\[
\lambda_1(x,u,y,v)=K_1-\frac{||gradf_2||^2}{f_1^2}.
\]
\end{thm}
\begin{proof}
Since $({}_{f_2}M_1\times{}_{f_1}M_2,F)$ is a space of scalar flag
curvature $\lambda_1(x,u,y,v)$, then we have
\be
\textbf{R}^{\ i}_{j\
kl}=\lambda_1(x,u,y,v)(\delta^i_lg_{jk}-\delta^i_{k}g_{jl}).\label{R1}
\ee
By setting (\ref{R1}) in (\ref{con1}), we obtain
\be
R^{\ i}_{j\ kl}=\Big[\lambda_1(x,u,y,v)+\frac{||gradf_2||^2}{f_1^2}\Big](\delta^i_lg_{jk}-\delta^i_{k}g_{jl}).\label{R2}
\ee
By using (\ref{R2}), the proof is completes.
\end{proof}

Similarly, we have the following.

\begin{thm}
Let $(M_2,F_2)$ be a Riemannian manifold and
$({}_{f_2}M_1\times{}_{f_1}M_2,F)$ be a doubly warped product
Finsler space of scalar flag curvature $\lambda_2(x,u,y,v)$. Then
$(M_2, F_2)$ has constant curvature $K_2$ if and only if
\[
\lambda_2(x,u,y,v)=K_2-\frac{||gradf_1||^2}{f_2^2}.
\]
\end{thm}

\smallskip

\begin{cor}
Let $({}_{f_2}M_1\times{}_{f_1}M_2,F)$ be a doubly warped product
Riemannian manifold of the constant curvature $\lambda$. Then
\begin{description}
  \item[(i)] if $f_1$ is constant on $M_1$, then $M_1$ and $M_2$ have constant curvatures $\lambda+||gradf_2||^2$ and
$\lambda$, respectively;
  \item[(ii)] if $f_2$ is constant on $M_2$, then $M_1$ and $M_2$ have constant curvatures $\lambda$ and
$\lambda+||gradf_1||^2$, respectively.
\end{description}
\end{cor}
\section{Doubly Warped Sasaki-Matsumoto Metric}
Let $(M, F)$ be a Finsler manifold. It is well known that there are
several ways to associate the slit tangent bundle $TM^\circ$ of $M$ with Riemannian
metrics which are naturally induced by the Finsler metric $F$. The most well-known such metric is the Sasaki-Matsumoto lift
\[
G=g_{ij}dx^i\otimes dx^j+g_{ij}\delta y^i\otimes\delta y^j.
\]
to the $TM^\circ$ (see \cite{PT1, PT2, PT3, PTZ}). Now, let $({}_{f_1}M_1\times{}_{f_2}M_2,F)$ be a doubly warped product Finsler manifold. Then the doubly warped Sasaki-Matsumoto metric can introduced as follows
\begin{eqnarray}
\nonumber \textbf{G}=f_2^2g_{ij}dx^i\otimes dx^j+f_1^2g_{\alpha\beta}du^\alpha\otimes du^\beta\!\!\!\!&+&\!\!\!\!\ f_2^2g_{ij}\delta^d y^i\otimes\delta^d y^j\\
\!\!\!\!&+&\!\!\!\!\ f_1^2g_{\alpha\beta}\delta^d v^\alpha\otimes\delta^d v^\beta,\label{Doubly}
\end{eqnarray}
where $\delta^d y^i$ and $\delta^d v^\alpha$ are defined by (\ref{new}).
\begin{prop}
Let $({}_{f_2}M_1\times{}_{f_1}M_2,F)$ be a doubly warped product
Finsler manifold. Then the Levi-Civita connection $\nabla^d$ on
the Riemannian manifold $(TM^\circ,G^d)$ is locally expressed as
follows:
\begin{eqnarray}
\nabla^d_{\frac{\D}{\D x^i}}\frac{\D}{\D x^j}\!\!\!\!&=&\!\!\!\!\ \textbf{F}^s_{ij}\frac{\D}{\D x^s}+(\frac{1}{2}\textbf{R}^s_{ij}-C^s_{ij})\frac{\partial}{\partial y^s}+\textbf{F}^\gamma_{ij}\frac{\D}{\D
u^\gamma}+\frac{1}{2}\textbf{R}^\gamma_{\ ij}\frac{\partial}{\partial
v^\gamma}\label{Levi1}\\
\nabla^d_{\frac{\D}{\D x^i}}\frac{\partial}{\partial y^j}\!\!\!\!&=&\!\!\!\!(C^s_{ij}+\frac{g_{rj}}{2}g^{ks}\textbf{R}^r_{ki})\frac{\D}{\D x^s}+\frac{f_2^2}{2f_1^2}g_{rj}g^{\gamma\mu}\textbf{R}^r_{\ \mu i}\frac{\D}{\D u^\gamma}
\nonumber\\
\!\!\!\!&&\!\!\!\!+\frac{g^{ks}}{2}(\frac{\D}{\D x^i}g_{jk}+\textbf{G}^r_{ij}g_{rk}-\textbf{G}^r_{ik}g_{rj})\frac{\partial}{\partial y^s}\nonumber\\
\!\!\!\!&&\!\!\!\!
+\frac{1}{2f_1^2}g^{\gamma\mu}(f_1^2\textbf{G}^\lambda_{ij}g_{\lambda\mu}-f_2^2\textbf{G}^r_{i\mu}g_{rj})
\frac{\partial}{\partial
v^\gamma}\label{Levi2}\\
\nabla^d_{\frac{\D}{\D x^i}}\frac{\D}{\D u^\beta}\!\!\!\!&=&\!\!\!\!\textbf{F}^s_{i\beta}\frac{\D}{\D x^s}+\frac{1}{2}\textbf{R}^s_{\ i\beta}\frac{\partial}{\partial y^s}+\textbf{F}^\gamma_{i\beta}\frac{\D}{\D
u^\gamma}+\frac{1}{2}\textbf{R}^\gamma_{\ i\beta}\frac{\partial}{\partial
v^\gamma}\label{Levi3}\\
\nabla^d_{\frac{\D}{\D x^i}}\frac{\partial}{\partial v^\beta}\!\!\!\!&=&\!\!\!\!\frac{f_1^2}{2f_2^2}g_{\lambda\beta}g^{ks}\textbf{R}^\lambda_{\ ki}\frac{\D}{\D x^s}+\frac{1}{2f_2^2}g^{ks}(f_2^2\textbf{G}^r_{i\beta}g_{rk}\nonumber\\
\!\!\!\!&&\!\!\!\!-f_1^2\textbf{G}^\lambda_{ik}g_{\lambda\beta})\frac{\partial}{\partial y^s}+\frac{1}{2}g_{\lambda\beta}g^{\mu\gamma}\textbf{R}^\lambda_{\ \mu i}\frac{\D}{\D u^\gamma}\nonumber\\
\!\!\!\!&&\!\!\!\! +\frac{1}{2f_1^2}g^{\gamma\mu}(\frac{\D}{\D
x^i}f_1^2g_{\beta\mu}+f_1^2\textbf{G}^\lambda_{i\beta}g_{\lambda\mu}
-f_1^2\textbf{G}^\lambda_{i\mu}g_{\lambda\beta})\frac{\partial}{\partial
v^\gamma}\label{Levi4}\\
\nabla^d_{\frac{\D}{\D u^\alpha}}\frac{\D}{\D u^\beta}\!\!\!\!&=&\!\!\!\!\textbf{F}^s_{\alpha\beta}\frac{\D}{\D x^s}+\frac{\textbf{R}^s_{\alpha\beta}}{2}\frac{\partial}{\partial y^s}+\textbf{F}^\gamma_{\alpha\beta}\frac{\D}{\D
u^\gamma}+(\frac{\textbf{R}^\gamma_{\alpha\beta}}{2}-C^\gamma_{\alpha\beta})\frac{\partial}{\partial
v^\gamma}\label{Levi5}\\
\nabla^d_{\frac{\D}{\D u^\alpha}}\frac{\D}{\D x^j}\!\!\!\!&=&\!\!\!\!\textbf{F}^s_{\alpha j}\frac{\D}{\D x^s}+\frac{1}{2}\textbf{R}^s_{\alpha j}\frac{\partial}{\partial y^s}+\textbf{F}^\gamma_{\alpha j}\frac{\D}{\D
u^\gamma}+\frac{1}{2}\textbf{R}^\gamma_{\alpha
j}\frac{\partial}{\partial v^\gamma}\label{Levi6}
\end{eqnarray}
\begin{eqnarray}
\nabla^d_{\frac{\D}{\D u^\alpha}}\frac{\partial}{\partial y^j}\!\!\!\!&=&\!\!\!\!\frac{1}{2f_2^2}g^{ks}(\frac{\D}{\D u^\alpha}f_2^2g_{jk}+f_2^2\textbf{G}^r_{\alpha j}g_{rk}-f_2^2\textbf{G}^r_{\alpha k}g_{rj})\frac{\partial}{\partial y^s}\nonumber\\
\!\!\!\!&&\!\!\!\!+\frac{1}{2}g_{rj}g^{ks}\textbf{R}^r_{k\alpha}\frac{\D}{\D x^s}+\frac{f_2^2}{2f_1^2}g_{rj}g^{\gamma\mu}\textbf{R}^r_{\mu \alpha}\frac{\D}{\D u^\gamma}\nonumber\\
\!\!\!\!&&\!\!\!\!
+\frac{1}{2f_1^2}g^{\gamma\mu}(f_1^2\textbf{G}^\lambda_{\alpha
j}g_{\lambda\mu}-f_2^2\textbf{G}^r_{\alpha\mu}g_{rj})\frac{\partial}{\partial
v^\gamma}\label{Levi7}\\
\nabla^d_{\frac{\D}{\D u^\alpha}}\frac{\partial}{\partial v^\beta}\!\!\!\!&=&\!\!\!\!\frac{f_1^2}{2f_2^2}g_{\lambda\beta}g^{ks}\textbf{R}^\lambda_{k\alpha}\frac{\D}{\D x^s}+\frac{1}{2f_2^2}g^{ks}(f_2^2\textbf{G}^r_{\alpha\beta}g_{rk}\nonumber\\
\!\!\!\!&&\!\!\!\!-f_1^2\textbf{G}^\lambda_{\alpha k}g_{\lambda\beta})\frac{\partial}{\partial y^s}+(C^\gamma_{\alpha\beta}+\frac{1}{2}g_{\lambda\beta}g^{\mu\gamma}\textbf{R}^\lambda_{\mu \alpha})\frac{\D}{\D u^\gamma}\nonumber\\
\!\!\!\!&&\!\!\!\! +\frac{1}{2f_1^2}g^{\gamma\mu}(\frac{\D}{\D
u^\alpha}f_1^2g_{\beta\mu}+f_1^2\textbf{G}^\lambda_{\alpha\beta}g_{\lambda\mu}-f_1^2\textbf{G}^\lambda_{\alpha\mu}
g_{\lambda\beta})\frac{\partial}{\partial v^\gamma}\label{Levi8}
\end{eqnarray}
\begin{eqnarray}
\nabla^d_{\frac{\partial}{\partial y^i }}\frac{\partial}{\partial y^j}\!\!\!\!&=&\!\!\!\!\frac{g^{ks}}{2}(
\textbf{G}^r_{kj}g_{ri}+\textbf{G}^r_{ki}g_{rj}-\frac{\D}{\D x^k}g_{ij})\frac{\D}{\D x^s}+C^s_{ij}\frac{\partial}{\partial y^s}\nonumber\\
\!\!\!\!&&\!\!\!\!
+\frac{1}{2f_1^2}g^{\gamma\mu}(f_2^2\textbf{G}^r_{\mu
j}g_{ri}+f_2^2\textbf{G}^r_{\mu i}g_{rj}-\frac{\D}{\D u^\mu}f_2^2g_{ij})\frac{\D}{\D
u^\gamma}\label{Levi9}\\
\nabla^d_{\frac{\partial}{\partial v^\alpha}}\frac{\partial}{\partial y^j}\!\!\!\!&=&\!\!\!\!\frac{1}{2f_2^2}g^{ks}(f_1^2\textbf{G}^\lambda_{k j}g_{\lambda\alpha}+f_2^2\textbf{G}^r_{k\alpha}g_{rj})\frac{\D}{\D x^s}\nonumber\\
\!\!\!\!&&\!\!\!\!
+\frac{1}{2f_1^2}g^{\gamma\mu}(f_1^2\textbf{G}^\lambda_{\mu
j}g_{\lambda\alpha}+f_2^2\textbf{G}^r_{\mu\alpha}g_{rj})\frac{\D}{\D
u^\gamma}=\nabla^d_{\frac{\partial}{\partial y^j}}\frac{\partial}{\partial v^\alpha}\label{Levi10}\\
\nabla^d_{\frac{\partial}{\partial v^\alpha}}\frac{\partial}{\partial v^\beta}\!\!\!\!&=&\!\!\!\!\frac{1}{2f_2^2}g^{ks}(-\frac{\D}{\D x^k}f_1^2g_{\alpha\beta}+f_1^2\textbf{G}^\lambda_{k \beta}g_{\lambda\alpha}+f_1^2\textbf{G}^\lambda_{k\alpha}g_{\lambda\beta})\frac{\D}{\D x^s}\nonumber\\
\!\!\!\!&&\!\!\!\!\
+\frac{1}{2f_1^2}g^{\gamma\mu}(-f_1^2\frac{\D}{\D u^\mu}g_{\alpha\beta}+f_1^2\textbf{G}^\lambda_{\mu
\beta}g_{\lambda\alpha}
+f_1^2\textbf{G}^\lambda_{\mu\alpha}g_{\lambda\beta})\frac{\D}{\D
u^\gamma}\nonumber\\
\!\!\!\!&&\!\!\!\!\ +C^\gamma_{\alpha\beta}\frac{\partial}{\partial
v^\gamma}\label{Levi11}
\end{eqnarray}
\end{prop}
The Levi-Civita connection $\nabla^d$ induces a connection
$\nabla$ on $VTM^\circ$, i.e.,
\begin{equation}
\nabla_Xv^dY=v^d(\nabla^d_Xv^dY),\label{ver}
\end{equation}
for $X,Y\in\Gamma(TTM^\circ)$. Then we have
\begin{lem}\label{Civita}
we have
\begin{eqnarray}
\nabla_{\frac{\D}{\D x^i}}\frac{\partial}{\partial
y^j}=\textbf{F}^s_{ij}\frac{\partial}{\partial
y^s}+\textbf{F}^\gamma_{ij}\frac{\partial}{\partial v^\gamma},\ \
\ \nabla_{\frac{\D}{\D u^\alpha}}\frac{\partial}{\partial
y^j}=\textbf{F}^s_{\alpha j}\frac{\partial}{\partial
y^s}+\textbf{F}^\gamma_{\alpha j}\frac{\partial}{\partial
v^\gamma}\label{mohem}
\\
\nabla_{\frac{\D}{\D x^i}}\frac{\partial}{\partial
v^\beta}=\textbf{F}^s_{i\beta}\frac{\partial}{\partial
y^s}+\textbf{F}^\gamma_{i\beta}\frac{\partial}{\partial
v^\gamma},\ \  \nabla_{\frac{\D}{\D
u^\alpha}}\frac{\partial}{\partial v^\beta}=\textbf{F}^s_{\alpha
\beta}\frac{\partial}{\partial y^s}+\textbf{F}^\gamma_{\alpha
\beta}\frac{\partial}{\partial v^\gamma}
\\
\nabla^d_{\frac{\partial}{\partial y^i}}\frac{\partial}{\partial
y^j}=C^s_{ij}\frac{\partial}{\partial y^s},\
\nabla^d_{\frac{\partial}{\partial
v^\alpha}}\frac{\partial}{\partial
v^\beta}=C^\gamma_{\alpha\beta}\frac{\partial}{\partial
v^\gamma},\ \nabla^d_{\frac{\partial}{\partial
y^i}}\frac{\partial}{\partial
v^\beta}=\nabla^d_{\frac{\partial}{\partial
v^\alpha}}\frac{\partial}{\partial y^j}=0.
\end{eqnarray}
\end{lem}
\begin{proof}
By using (\ref{Levi2}), (\ref{Levi4}),
(\ref{Levi7})-(\ref{Levi11}) and (\ref{ver}) it is sufficient to
proof the following equation:
\begin{equation}
\frac{1}{2}\textbf{g}^{ce}(\frac{\D \textbf{g}_{ea}}{\D
\textbf{x}^b}+\textbf{g}_{de}\textbf{G}^d_{ba}-\textbf{g}_{da}\textbf{G}^d_{be})=\textbf{F}^c_{ab}.\label{pr}
\end{equation}
By using (i) of (\ref{spray6}), we derive
\begin{eqnarray}
\textbf{g}_{de}\textbf{G}^d_{ba}-\textbf{g}_{da}\textbf{G}^d_{be}\!\!\!\!&=&\!\!\!\!\ \textbf{g}_{de}\frac{\partial^2\textbf{G}^d}{\partial
\textbf{y}^a\partial \textbf{y}^b}-\textbf{g}_{da}\frac{\partial^2\textbf{G}^d}{\partial
\textbf{y}^e\partial \textbf{y}^b}\nonumber\\
\!\!\!\!&=&\!\!\!\!\ \frac{\partial}{\partial \textbf{y}^b}(\frac{\partial
\textbf{G}_e}{\partial \textbf{y}^a}-\frac{\partial
\textbf{G}_a}{\partial \textbf{y}^e})-\textbf{G}^d_a\frac{\partial
\textbf{g}_{de}}{\partial
\textbf{y}^b}+\textbf{G}^d_e\frac{\partial
\textbf{g}_{da}}{\partial \textbf{y}^b},\label{pr1}
\end{eqnarray}
where $\textbf{G}_e:=\textbf{g}_{de}\textbf{G}^d$. By direct
calculations and using (\ref{Mat}) and (\ref{spray3}), we deduce
that
\begin{equation}
\frac{\partial}{\partial \textbf{y}^b}(\frac{\partial
\textbf{G}_e}{\partial \textbf{y}^a}-\frac{\partial
\textbf{G}_a}{\partial \textbf{y}^e})=\frac{\partial
\textbf{g}_{eb}}{\partial \textbf{x}^a}-\frac{\partial
\textbf{g}_{ab}}{\partial \textbf{x}^e}.\label{pr2}
\end{equation}
Putting (\ref{pr2}) in (\ref{pr1}) implies that
\begin{equation}
\textbf{g}_{de}\textbf{G}^d_{ba}-\textbf{g}_{da}\textbf{G}^d_{be}=\frac{\D
\textbf{g}_{eb}}{\D \textbf{x}^a}-\frac{\D \textbf{g}_{ab}}{\D
\textbf{x}^e}.\label{g1}
\end{equation}
By plugging (\ref{g1}) in the left side of equation (\ref{pr})
and using (\ref{hor}),  we obtain the right side of equation
(\ref{pr}).
\end{proof}
We say that the vertical distribution $VTM^\circ$ is totally
geodesic in $TTM^\circ$ if
\[
\nabla^d_{\frac{\partial}{\partial \textbf{y}^a}}\frac{\partial}{\partial \textbf{y}^b}\in\Gamma(VTM^\circ).
\]
Similarly, we say that the horizontal distribution $HTM^\circ$ is totally geodesic in
$TTM^\circ$ if
\[
\nabla^d_{\frac{\D}{\D \textbf{x}^a}}\frac{\D}{\D
\textbf{x}^b}\in\Gamma(HTM^\circ).
\]
For more details, see \cite{wu}.
\begin{prop}
Let $({}_{f_2}M_1\times{}_{f_1}M_2,F)$ be a doubly warped product
Finsler manifold. Then $VTM^\circ$ is totally geodesic if and only
if $\textbf{F}^c_{ab}=\textbf{G}^c_{ab}$.
\end{prop}
\begin{proof}
By using the definition of totally geodesic, we deduce that
$VTM^\circ$ is totally geodesic if and only if
\begin{equation}
\nabla^d_{\frac{\partial}{\partial y^i}}\frac{\partial}{\partial
y^j}\in\Gamma(VTM^\circ),\ \ \ \nabla^d_{\frac{\partial}{\partial
y^i}}\frac{\partial}{\partial v^\beta}\in\Gamma(VTM^\circ),
\end{equation}
\begin{equation}
\nabla^d_{\frac{\partial}{\partial
v^\alpha}}\frac{\partial}{\partial y^j}\in\Gamma(VTM^\circ),\ \ \
\nabla^d_{\frac{\partial}{\partial
v^\alpha}}\frac{\partial}{\partial v^\beta}\in\Gamma(VTM^\circ).
\end{equation}
Since $G^d$ is parallel with respect to $\nabla^d$, then we have
\begin{equation}
\textbf{G}(\nabla^d_{\frac{\partial}{\partial
y^i}}\frac{\partial}{\partial y^j}, \frac{\D}{\D
x^h})+\textbf{G}(\nabla^d_{\frac{\partial}{\partial y^i}}\frac{\D}{\D
x^h},\frac{\partial}{\partial y^j})=0,
\end{equation}
\begin{equation}
\textbf{G}(\nabla^d_{\frac{\partial}{\partial
y^i}}\frac{\partial}{\partial y^j}, \frac{\D}{\D
u^\lambda})+\textbf{G}(\nabla^d_{\frac{\partial}{\partial
y^i}}\frac{\D}{\D u^\lambda},\frac{\partial}{\partial y^j})=0.
\end{equation}
By using (\ref{Levi9}) and (\ref{mohem}), we obtain
\begin{equation}
\frac{1}{2}g^{kl}(-\frac{\D}{\D x^k}g_{ij}
+\textbf{G}^r_{kj}g_{ri}+\textbf{G}^r_{ki}g_{rj})=g^{lh}(\textbf{G}^s_{ih}-\textbf{F}^s_{ih})g_{sj},
\end{equation}
\begin{equation}
\frac{1}{2}g^{\nu\mu}(-\frac{\D}{\D
u^\mu}f_2^2g_{ij}+f_2^2\textbf{G}^r_{\mu
j}g_{ri}+f_2^2\textbf{G}^r_{\mu
i}g_{rj})=f_2^2g^{\lambda\nu}(\textbf{G}^s_{i\lambda}-\textbf{F}^s_{i\lambda})g_{sj}.
\end{equation}
Putting the above equations in (\ref{Levi9}) give us
\[
\nabla^d_{\frac{\partial}{\partial y^i }}\frac{\partial}{\partial
y^j}=[g^{sh}(\textbf{G}^l_{ih}-\textbf{F}^l_{ih})\frac{\D}{\D
x^s}+\frac{f_2^2}{f_1^2}g^{\lambda\gamma}(\textbf{G}^l_{i\lambda}-\textbf{F}^l_{i\lambda})\frac{\D}{\D
u^\gamma}]g_{lj}+C^s_{ij}\frac{\partial}{\partial
y^s}.
\]
Therefore $\nabla^d_{\frac{\partial}{\partial y^i}}\frac{\partial}{\partial y^j}\in\Gamma(VTM^\circ)$ if and
only if $\textbf{G}^l_{ih}=\textbf{F}^l_{ih}$ and $\textbf{G}^l_{i\lambda}=\textbf{F}^l_{i\lambda}$. Similarly,  we
obtain the another relations.
\end{proof}
\begin{cor}
Let $({}_{f_2}M_1\times{}_{f_1}M_2,F)$ be a doubly warped product
Riemannian manifold. Then $VTM^\circ$ is totally geodesic distribution.
\end{cor}

\bigskip

\begin{prop}
Let $({}_{f_2}M_1\times{}_{f_1}M_2,F)$ be a doubly warped product
Finsler manifold. Then $HTM^\circ$ is totally geodesic if and only
if $(M_1, F_1)$ and $(M_2,F_2)$ are Riemannian manifolds and
$\textbf{R}^\gamma_{\ ij}=\textbf{R}^s_{\
i\beta}=\textbf{R}^\gamma_{\ i\beta}=\textbf{R}^s_{\
\alpha\beta}=0$.
\end{prop}
\begin{proof}
By  definition, $HTM^\circ$ is totally geodesic if and only if
\begin{equation}
\nabla^d_{\frac{\D}{\D x^i}}\frac{\D}{\D
x^j}\in\Gamma(HTM^\circ),\ \ \ \nabla^d_{\frac{\D}{\D
x^i}}\frac{\D}{\D u^\beta}\in\Gamma(HTM^\circ),
\end{equation}
\begin{equation}
\nabla^d_{\frac{\D}{\D u^\alpha}}\frac{\D}{\D
x^j}\in\Gamma(HTM^\circ),\ \ \ \nabla^d_{\frac{\D}{\D
u^\alpha}}\frac{\D}{\D u^\beta}\in\Gamma(HTM^\circ).
\end{equation}
(\ref{Levi1}) implies that $\nabla^d_{\frac{\D}{\D
x^i}}\frac{\D}{\D x^j}\in\Gamma(HTM^\circ)$ if and only if
$\textbf{R}^\gamma_{\ ij}=0$ and
\begin{equation}
-C^s_{ij}+\frac{1}{2}\textbf{R}^s_{\ ij}=0.\label{car1}
\end{equation}
Interchanging $i$ and $j$ in the above equation gives us
\begin{equation}
-C^s_{ji}+\frac{1}{2}\textbf{R}^s_{\ ji}=0.\label{car2}
\end{equation}
It is remarkable that  $C^s_{ij}$ and $\textbf{R}^s_{\ ij}$ are symmetric and
skew-symmetric tensors with respect to $i$ and $j$, respectively. Then (\ref{car1})+(\ref{car2}), implies that $C^s_{ij}=0$, i.e., $(M_1, F_1)$ is a Riemannian manifold. In the similar way, we can prove another relations.
\end{proof}
%---------------------------------------------------------------------------------------------------------
\section{Doubly Warped Vaisman Connection}
%---------------------------------------------------------------------------------------------------------
In this section, the Riemannian manifold $(TM^0,\textbf{G})$ is considered for which $M={}_{f_2}M_1\times{}_{f_1}M_2$ and $\textbf{G}$ is given by (\ref{Doubly}) and the vertical foliation $\mathcal{F}_V$ (i.e., $TTM^\circ=HTM^\circ\oplus VTM^\circ$) on it.
Also, we consider the notation from \cite{Balan} and \cite{Vais},
related to foliated manifolds entitled \textit{Vaisman
connection}. The Vaisman connection $\nabla^v$ on the Riemannian
foliated manifold $(TM^\circ, \mathcal{F}_V, \textbf{G})$, is
uniquely defined by the following conditions:
\begin{description}
\item[(i)] if $Y\in\Gamma(VTM^\circ)$ (respectively
$\in\Gamma(HTM^\circ)$), then $\nabla^v_XY\in\Gamma(VTM^\circ)$
(respectively $\in\Gamma(HTM^\circ)$) for every $X$;
\item[(ii)] if $X, Y, Z\in\Gamma(VTM^\circ)$
($\Gamma(HTM^\circ)$), then $(\nabla^v_XG)(Y,Z)=0$;
\item[(iii)] $v^d(T(X,Y))=0$ if at least one of the arguments is in
$\Gamma(VTM^\circ)$ and $h^d(T(X,Y))=0$ if at least one of the
arguments is in $\Gamma(HTM^\circ)$.
\end{description}

Now, we are going to compute the  Vaisman connection  on the Riemannian
foliated manifold $(TM^\circ, \mathcal{F}_V, \textbf{G})$.

\begin{prop}\label{Vaisman}
Let $({}_{f_2}M_1\times{}_{f_1}M_2,F)$ be a doubly warped product
Finsler manifold. Then the Vaisman connection $\nabla^v$ on
$(TM^\circ, \mathcal{F}_V, \textbf{G})$, is locally expressed
with respect to the adapted local basis $\{\frac{\D}{\D x^i},
\frac{\D}{\D u^\alpha}, \frac{\partial}{\partial y^i},
\frac{\partial}{\partial v^\alpha}\}$ as follows:
\begin{eqnarray}
\nabla^v_{\frac{\D}{\D x^i}}\frac{\partial}{\partial
y^j}=\textbf{G}^k_{ij}\frac{\partial}{\partial
y^k}+\textbf{G}^\gamma_{ij}\frac{\partial}{\partial v^\gamma}, \
\nabla^v_{\frac{\D}{\D u^\alpha}}\frac{\partial}{\partial v^\beta}
=\textbf{G}^k_{\alpha\beta}\frac{\partial}{\partial
y^k}+\textbf{G}^\gamma_{\alpha\beta}\frac{\partial}{\partial
v^\gamma}\label{Vais1}
\\
\nabla^v_{\frac{\D}{\D u^\alpha}}\frac{\partial}{\partial
y^j}=\textbf{G}^k_{\alpha j}\frac{\partial}{\partial
y^k}+\textbf{G}^\gamma_{\alpha j}\frac{\partial}{\partial
v^\gamma}, \ \nabla^v_{\frac{\D}{\D x^i}}\frac{\partial}{\partial
v^\beta} =\textbf{G}^k_{i\beta}\frac{\partial}{\partial
y^k}+\textbf{G}^\gamma_{i\beta}\frac{\partial}{\partial
v^\gamma}\label{Vais2}
\end{eqnarray}
\begin{equation}
\nabla^v_{\frac{\D}{\D x^i}}\frac{\D}{\D x^j}=\textbf{F}^k_{
ij}\frac{\D}{\D x^k}+\textbf{F}^\gamma_{ij}\frac{\D}{\D u^\gamma},
\ \nabla^v_{\frac{\D}{\D u^\alpha}}\frac{\D}{\D
u^\beta}=\textbf{F}^k_{ \alpha\beta}\frac{\D}{\D
x^k}+\textbf{F}^\gamma_{\alpha\beta}\frac{\D}{\D
u^\gamma},\label{Vais3}
\end{equation}
\begin{equation}
\nabla^v_{\frac{\D}{\D x^i}}\frac{\D}{\D u^\beta}=\textbf{F}^k_{
i\beta}\frac{\D}{\D x^k}+\textbf{F}^\gamma_{i\beta}\frac{\D}{\D
u^\gamma}, \ \nabla^v_{\frac{\D}{\D u^\alpha}}\frac{\D}{\D
x^j}=\textbf{F}^k_{ \alpha j}\frac{\D}{\D
x^k}+\textbf{F}^\gamma_{\alpha j}\frac{\D}{\D
u^\gamma},\label{Vais4}
\end{equation}
\begin{equation}
\nabla^v_{\frac{\partial}{\partial y^i}}\frac{\partial}{\partial
y^j} =C^k_{ij}\frac{\partial}{\partial y^k},\
\nabla^v_{\frac{\partial}{\partial
v^\alpha}}\frac{\partial}{\partial
v^\beta}=C^\gamma_{\alpha\beta}\frac{\partial}{\partial
v^\gamma},\ \nabla^v_{\frac{\partial}{\partial
y^i}}\frac{\partial}{\partial v^\beta}=
\nabla^v_{\frac{\partial}{\partial
v^\alpha}}\frac{\partial}{\partial y^j}=0,\label{Vais5}
\end{equation}
\begin{equation}
\nabla^v_{\frac{\partial}{\partial y^i}}\frac{\D}{\D
x^j}=\nabla^v_{\frac{\partial}{\partial v^\alpha}}\frac{\D}{\D
x^j}=\nabla^v_{\frac{\partial}{\partial v^\alpha}}\frac{\D}{\D
u^\beta}=\nabla^v_{\frac{\partial}{\partial y^i}}\frac{\D}{\D
u^\beta}=0\label{Vais6}.
\end{equation}
\end{prop}
\begin{proof}
From the condition (i) of Vaisman connection, we have
\be
\nabla^v_{\frac{\D}{\D x^i}}\frac{\partial}{\partial
y^j}=A^k_{ij}\frac{\partial}{\partial
y^k}+A^\gamma_{ij}\frac{\partial}{\partial v^\gamma},\ \ \
\nabla^v_{\frac{\partial}{\partial y^j}}\frac{\D}{\D
x^i}=B^k_{ji}\frac{\D}{\D x^k}+B^\gamma_{ji}\frac{\D}{\D
u^\gamma}.\label{V}
\ee
By using (\ref{V}) and the condition (ii) of Vaisman connection, we get
\[
0=V(T(\frac{\D}{\D x^i}, \frac{\partial}{\partial
y^j}))=(A^k_{ij}-\textbf{G}^k_{ij})\frac{\partial}{\partial
y^k}+(A^\gamma_{ij}-\textbf{G}^\gamma_{ij})\frac{\partial}{\partial
v^\gamma}-B^k_{ji}\frac{\D}{\D x^k}-B^\gamma_{ji}\frac{\D}{\D
u^\gamma}.
\]
The above equation $A^k_{ij}=\textbf{G}^k_{ij}$,
$A^\gamma_{ij}=\textbf{G}^\gamma_{ij}$ and
$B^k_{ji}=B^\gamma_{ji}=0$. Therefore, we obtain the first equation
of (\ref{Vais1}) and the first equation of (\ref{Vais6}).
Similarly, we can obtain the another relations.
\end{proof}
The Lemma \ref{Civita} and Proposition \ref{Vaisman}, give us the following.
\begin{thm}
Let $({}_{f_2}M_1\times{}_{f_1}M_2,F)$ be a doubly warped product
Finsler manifold. Then the Levi-Civita and the Vaisman connections
on the foliated manifold $(TM^\circ, \mathcal{F}_V, \textbf{G})$ induce
the same connection on the structural bundle if and only if
$\textbf{F}^c_{ab}=\textbf{G}^c_{ab}$.
\end{thm}

Therefore, we conclude the following.

\begin{cor}
Let $({}_{f_2}M_1\times{}_{f_1}M_2,F)$ be a doubly warped product
Riemannian manifold. Then the Levi-Civita and the Vaisman
connections on the foliated manifold $(TM^\circ, \mathcal{F}_V,
\textbf{G})$ induce the same connection on the structural bundle.
\end{cor}

\smallskip

\begin{defn}
\emph{A Riemannian foliated manifold with the Riemannian metric $G$ is
called a \textit{Reinhart} space if and only if
\begin{equation}
(\nabla^v_XG)(Y,Z)=0,\label{Rein}
\end{equation}
for all the sections $X$ of the structural bundle and $Y, Z$
sections of the \textit{transversal bundle}, where the covariant
derivative is taken with respect to the Vaisman connection of the
manifold \cite{TP}.}
\end{defn}
\begin{thm}\label{TH4}
Let $({}_{f_1}M_1\times{}_{f_2}M_2,F)$ be a doubly warped product
Finsler manifold. The foliated manifold $(TM^\circ, \mathcal{F}_V,
\textbf{G})$ is a Reinhart space if and only if $(M_1,F_1)$ and
$(M_2,F_2)$ are Riemannian manifolds.
\end{thm}
\begin{proof}
Let $X=X^i\frac{\partial}{\partial
y^i}+X^\alpha\frac{\partial}{\partial
v^\alpha}\in\Gamma(VTM^\circ)$ and $Y=Y^j\frac{\D}{\D
x^i}+Y^\beta\frac{\D}{\D u^\beta}$, $Z=Z^k\frac{\D}{\D
x^k}+Z^\gamma\frac{\D}{\D u^\gamma}$ belong to
$\Gamma(HTM^\circ)$, also $\nabla^v$ be the Vaisman connection on
$(TM^\circ, \mathcal{F}_V, G^d)$. By  (\ref{Vais6}), we obtain
\begin{eqnarray*}
(\nabla^v_X\textbf{G})(Y,Z)\!\!\!\!&=&\!\!\!\!X^i\frac{\partial}{\partial
y^i}(Y^jZ^kf_2^2g_{jk})+X^i\frac{\partial}{\partial y^i}(Y^\beta
Z^\gamma
f_1^2g_{\beta\gamma})\nonumber\\
\!\!\!\!&&\!\!\!\!+X^\alpha\frac{\partial}{\partial v^\alpha}(Y^j
Z^k f_2^2g_{jk})+X^\alpha\frac{\partial}{\partial
v^\alpha}(Y^\beta Z^\gamma f_1^2g_{\beta\gamma})\nonumber\\
\!\!\!\!&&\!\!\!\!-X^i\frac{\partial Y^j}{\partial
y^i}Z^kf_2^2g_{jk}-X^i\frac{\partial Y^\beta}{\partial
y^i}Z^\gamma f_1^2g_{\beta\gamma}-X^\alpha\frac{\partial
Y^j}{\partial v^\alpha}Z^kf_2^2g_{jk}\nonumber\\
\!\!\!\!&&\!\!\!\!-X^\alpha\frac{\partial Y^\beta}{\partial
v^\alpha}Z^\gamma f_1^2g_{\beta\gamma}-X^iY^j\frac{\partial
Z^k}{\partial y^i}f_2^2g_{jk}-X^iY^\beta\frac{\partial
Z^\gamma}{\partial y^i} f_1^2g_{\beta\gamma}\nonumber\\
\!\!\!\!&&\!\!\!\!-X^\alpha Y^j\frac{\partial Z^k}{\partial
v^\alpha}f_2^2g_{jk}-X^\alpha Y^\beta\frac{\partial
Z^\gamma}{\partial v^\alpha}f_1^2g_{\beta\gamma}\nonumber\\
\!\!\!\!&=&\!\!\!\!2X^iY^jZ^kf_2^2C_{ijk}+2X^\alpha Y^\beta
Z^\gamma f_1^2C_{\alpha\beta\gamma}.
\end{eqnarray*}
Hence the condition (\ref{Rein}) for all $Y,
Z\in\Gamma(HTM^\circ)$ is equivalent to $C_{ijk}=0$ and
$C_{\alpha\beta\gamma}=0$, which is equal to that $(M_1,F_1)$ and
$(M_2,F_2)$ are Riemannian manifolds.
\end{proof}
\section{K\"{a}hlerian Structures}
In this section, we define an almost complex structure on the slit tangent bundle of a doubly warped product Finsler manifold and show that this structure together the doubly warped Sasaki-Matsumoto metric construct an almost Hermitian structure. Then we find a condition, under which this structure can be a K\"{a}hler Structure.

We consider the $\digamma{(TM^\circ)}$- linear mapping $\textbf{J}:\chi(TM^\circ)\rightarrow \chi(TM^\circ)$, defined by
\[
\textbf{J}=\frac{\D}{\D x^i}\otimes\delta^d y^i-\frac{\partial}{\partial y^i}\otimes dx^i+\frac{\D}{\D u^\alpha}\otimes\delta^dv^\alpha-\frac{\partial}{\partial v^\alpha}\otimes du^\alpha,
\]
or
\begin{eqnarray*}
&&\textbf{J}(\frac{\D}{\D x^i})=-\frac{\partial}{\partial y^i},\ \ \ \  \textbf{J}(\frac{\partial}{\partial y^i})=\frac{\D}{\D x^i},\\ &&\textbf{J}(\frac{\D}{\D u^\alpha})=-\frac{\partial}{\partial v^\alpha},\ \ \
\textbf{J}(\frac{\partial}{\partial v^\alpha})=\frac{\D}{\D u^\alpha}.\label{complex}
\end{eqnarray*}
It is easy to see that $\textbf{J}^2=-I$, i.e., $\textbf{J}$ is an almost complex structure on $TM^\circ$. Also, simple calculations give us $\textbf{G}(\textbf{J}X, \textbf{J}Y)=\textbf{G}(X,Y)$, where $X, Y\in \Gamma(TM^\circ)$. It means that $\textbf{G}$ is almost Hermitian with respect to $\textbf{J}$. The almost symplectic structure associated to the almost Hermitian structure $(\textbf{G}, \textbf{J})$ is defined by
\[
\Omega(X,Y):=\textbf{G}(X,\textbf{J}Y),\, \, \, \, \forall X,Y\in \Gamma (TM^\circ).
\]
By using (\ref{complex}) and the above equation, we obtain
\begin{eqnarray*}
\Omega(\frac{\D}{\D x^i}, \frac{\partial}{\partial y^j}) \!\!\!\!&=&\!\!\!\! \textbf{G}(\frac{\D}{\D x^i}, \textbf{J}(\frac{\partial}{\partial y^j}))\\
\\ \!\!\!\!&=&\!\!\!\!\ \textbf{G}(\frac{\D}{\D x^i}, \frac{\D}{\D x^j})\\
\!\!\!\!&=&\!\!\!\!\ f_2^2g_{ij}.
\end{eqnarray*}
Similarly, we  get the following
\begin{eqnarray*}
\Omega(\frac{\D}{\D u^\alpha}, \frac{\partial}{\partial v^\beta})= f_1^2g_{\alpha\beta}
\end{eqnarray*}
and
\begin{eqnarray*}
\Omega(\frac{\D}{\D x^i}, \frac{\D}{\D x^j}) \!\!\!\!&=&\!\!\!\!\  \Omega(\frac{\D}{\D x^i}, \frac{\D}{\D u^\alpha})=\Omega(\frac{\D}{\D x^i}, \frac{\partial}{\partial v^\beta})\\
 \!\!\!\!&=&\!\!\!\!\ \Omega(\frac{\D}{\D u^\alpha}, \frac{\D}{\D u^\beta})=\Omega(\frac{\D}{\D u^\alpha}, \frac{\partial}{\partial v^\beta})\\
 \!\!\!\!&=&\!\!\!\!\ \Omega(\frac{\partial}{\partial y^i}, \frac{\partial}{\partial y^j})= \Omega(\frac{\partial}{\partial y^i}, \frac{\partial}{\partial v^\beta})\\
 \!\!\!\!&=&\!\!\!\!\ \Omega(\frac{\partial}{\partial v^\alpha}, \frac{\partial}{\partial v^\beta})=0.
\end{eqnarray*}
Therefore, we can rewrite $\Omega$ as follows:
\[
\Omega=f_2^2g_{ij}dx^i\wedge\delta^dy^j+f_1^2g_{\alpha\beta}du^\alpha\wedge\delta^dv^\beta.
\]
By direct calculations, it is result that $\Omega=d\omega$, where
\[
\omega=f_2^2y^ig_{ij}dx^j+f_1^2v^\alpha g_{\alpha\beta}du^\beta.
\]
Thus $\Omega$ is a close form. By attention to these explanations,  we can conclude the following theorem.
\begin{thm}\label{feri}
 $(TM^\circ,\textbf{G}, \textbf{J})$ is an almost K\"{a}hlerian manifold.
\end{thm}

\bigskip

Consequently, the K\"{a}hler structure on $(TM^\circ, \textbf{J}, \textbf{G})$ is equivalent to  the
integrability condition of $\textbf{J}$. The integrability of $\textbf{J}$ is equal to the vanishing
of tensor field $N_{\textbf{J}}$, which is given by following
\begin{equation}
N_\textbf{J}(X,Y)=[\textbf{J}X,\textbf{J}Y]-\textbf{J}[\textbf{J}X,Y]-\textbf{J}[X,\textbf{J}Y]-[X,Y],\label{N}
\end{equation}
where $X,Y \in\Gamma(TM^\circ)$. By (\ref{N}), in computing $N_\textbf{J}$, the following equations are
presented
\begin{eqnarray}
N_\textbf{J}(\frac{\D}{\D x^i}, \frac{\D}{\D x^j}) \!\!\!\!&=&\!\!\!\!\ -N_\textbf{J}(\frac{\partial}{\partial y^i}, \frac{\partial}{\partial y^j})=-\textbf{R}^k_{\ ij}\frac{\partial}{\partial y^k}-\textbf{R}^\gamma_{\ ij}\frac{\partial}{\partial v^\gamma}
\\
N_\textbf{J}(\frac{\D}{\D u^\alpha}, \frac{\D}{\D u^\beta}) \!\!\!\!&=&\!\!\!\!\ -N_\textbf{J}(\frac{\partial}{\partial v^\alpha}, \frac{\partial}{\partial v^\beta})=-\textbf{R}^k_{\ \alpha\beta}\frac{\partial}{\partial y^k}-\textbf{R}^\gamma_{\ \alpha\beta}\frac{\partial}{\partial v^\gamma}
\\
N_\textbf{J}(\frac{\D}{\D x^i},\frac{\partial}{\partial v^\alpha})\!\!\!\!&=&\!\!\!\!\ -N_\textbf{J}(\frac{\D}{\D u^\alpha}, \frac{\partial}{\partial y^i})=-\textbf{R}^k_{\ i\alpha}\frac{\D}{\D x^k}-\textbf{R}^\gamma_{\ i\alpha}\frac{\D}{\D u^\gamma}
\\
N_\textbf{J}(\frac{\D}{\D x^i}, \frac{\D}{\D u^\alpha}) \!\!\!\!&=&\!\!\!\!\ -N_\textbf{J}(\frac{\partial}{\partial v^\alpha}, \frac{\partial}{\partial y^i})=-\textbf{R}^k_{\ i\alpha}\frac{\partial}{\partial y^k}-\textbf{R}^\gamma_{\ i\alpha}\frac{\partial}{\partial v^\gamma}
\\
N_\textbf{J}(\frac{\D}{\D x^i},\frac{\partial}{\partial y^j}) \!\!\!\!&=&\!\!\!\!\ -\textbf{R}^k_{\ ij}\frac{\D}{\D x^k}-\textbf{R}^\gamma_{\ ij}\frac{\D}{\D u^\gamma}
\\
N_\textbf{J}(\frac{\D}{\D u^\alpha},\frac{\partial}{\partial v^\beta})\!\!\!\!&=&\!\!\!\!\ -\textbf{R}^k_{\ \alpha\beta}\frac{\D}{\D x^k}-\textbf{R}^\gamma_{\ \alpha\beta}\frac{\D}{\D u^\gamma}.
\end{eqnarray}
Thus we have the following.
\begin{lem}\label{esi}
The complex structure $\textbf{J}$ defined by (\ref{complex}) is integrable if and only if $\textbf{R}^a_{\ bc}=0$, where $a, b, c=1,\ldots,n_1+n_2$.
\end{lem}
On the other hand, $\textbf{R}^a_{\ bc}=0$ is equivalent to the integrability of $HTM^\circ$. Therefore, using Theorem \ref{feri} and Lemma \ref{esi}, we conclude the following.
\begin{thm}
$(TM^\circ,\textbf{G}, \textbf{J})$ is a K\"{a}hlerian manifold if and only if the doubly warped horizontal distribution $HTM^\circ$ is integrable.
\end{thm}
%--------------------------------------------------------------------------------------------------------------

%--------------------------------------------------------------------------------------------------------------

\noindent
Esmail Peyghan\\
Department of Mathematics, Faculty  of Science\\
Arak University\\
Arak 38156-8-8349,  Iran\\
Email: epeyghan@gmail.com

\bigskip

\noindent
Akbar Tayebi\\
Department of Mathematics, Faculty  of Science\\
Qom University\\
Qom, Iran\\
Email: akbar.tayebi@gmail.com

\end{document}